\documentclass[11pt]{article}      
\usepackage{graphicx}
\usepackage{color}
\usepackage{flafter}
\usepackage[centertags]{amsmath}
\usepackage{amsfonts}
\usepackage{amssymb}
\usepackage{cite}
\usepackage{enumerate}
\usepackage{tikz}
\usepackage{graphicx}
\usepackage{subfigure}
\usepackage{caption}

\usepackage[top=25.4mm, bottom=25.4mm, left=20.7mm, right=20.2mm]{geometry}
\allowdisplaybreaks
\begin{document}

\title{\textbf{A note on a conjecture of star chromatic index for outerplanar graphs}\thanks{This work is supported  by the Science \& Technology development Fund of Tianjin Education Commission for Higher Education, China (No. 2019KJ090).\newline{{\it $^a$Corresponding author. E-mail: dengyuqiu1980@126.com , xcdeng@mail.tjnu.edu.cn}}}}
 \author{Xingchao Deng$^a$,  Qingye Yao$^b$, Yanbing Zhang$^c$ and Xudong Cui$^d$\\
 \footnotesize $^a$ $^b$ $^c$ $^d$College of Mathematical Science, Tianjin  Normal University\\
 \footnotesize Tianjin City, 300071, P. R. China}
\date{}
\maketitle
%
\begin{abstract}

A {\it star edge coloring} of a graph $G$ is a proper edge coloring of $G$ without bichromatic paths or cycles of length four.
The {\it star chromatic index}, $\chi_{st}^{'} (G ),$ of $G$ is the minimum number $k$ for which $G$ has a star edge coloring by $k$ colors. In \cite{LB},
 L. Bezegov$\acute{a}$  et al.  conjectured that $\chi_{st}^{'} (G )\leq \lfloor\frac{3\Delta}{2}\rfloor+1$ when $G$ is an outerplanar graph with
 maximum degree $\Delta \geq 3.$ In this paper we obtained that $\chi_{st}^{'}(G) \leq \Delta+6$ when $G$ is an 2-connected  outerplanar graph with diameter 2 or 3. If $G$ is an 2-connected outerplanar graph with maximum degree 5, then $\chi_{st}^{'}(G) \leq 9.$
\end{abstract}
\section{Introduction}
All graphs in this paper are simple and undirected. For a graph $G,$ we use $V(G)$ and $E(G)$ to denote the vertex and edge set, $\left| V(G) \right|$ and $\left| E(G) \right|$ to denote the number of vertices and edges respectively. The {\it star coloring} of a graph is a proper coloring and the vertices of any two color classes of the star colored graph induce a star forest. In 1973,  Gr$\mathrm{\ddot{u}}$nbaum initiated the star coloring in \cite{BG}. A proper edge coloring of a graph $G,$ which has no 2-edge colored 4-length paths or cycles,  is a {\it star edge coloring} of $G.$  The minimum number of colors needed for a star edge coloring of $G$ is called the {\it star chromatic index} of it, denoted by $\chi^{\prime}_{st}(G)$.  In 2008, Xinsheng Liu and Kai Deng  introduced the star edge coloring  in \cite{XL}. Obviously, a star edge coloring of a graph $G$ is a star coloring of the line graph of $G$. A {\it strong edge coloring} of a graph $G,$ denoted by $\chi_{s}^{\prime}(G),$ is a proper edge coloring so that no edge can be adjacent to two edges with the same color in which every color class gives an induced matching. The definition was introduced by Fouquet and Jolivet (1983) to solve a problem involving radio networks and their frequencies. It is easy to know that the strong edge coloring of a graph $G$ is a star edge coloring of it, then $\chi^{'}_{st}(G)\leq\chi_{s}^{\prime}(G).$ The diameter of a graph $G,$ denoted by $diam(G),$ is the maximum distance of any two vertices in the graph $G.$  We use the symbol $G_1\vee G_2$ to denote the join graph of $G_1, G_2.$ A path with $n$ vertices is denoted by $P_n.$  A graph with $n$ vertices is denoted by $G_n.$ \\

If a graph $G$ can be drawn in the plane so that its edges intersect only at their ends, then $G$ is called a planar graph. In particular, when all vertices of a planar graph $G$ lie on a same face, then $G$ is an outerplanar graph. A minor of a graph $G$  is obtained by deleting some vertices, edges or contracting some edges of the graph $G.$ For the graphs $G$ and $H$ such that no minor of $G$ is isomorphic to $H,$  then $G$ is called a $H$ minor free graph.

Recently, Dvo$\mathrm{\check{r}}$$\mathrm{\acute{a}}$ and  Mohar have studied the star edge coloring of complete graphs and 3-regular graphs. They obtained some results and proposed a conjecture shown as follows:\\

\noindent{\bf Theorem 1.1\cite{ZD}}
Let $K_n$ be a complete graph with $n$ vertices, then
\begin{align*}
2n(1+o(1))\leq \chi_{st}^{'}(K_n)\leq n\frac{2^{2\sqrt{2}(1+o(1))\sqrt{\log n}}}{(\log n^{\frac{1}{4}})} O(\frac{\log\Delta}{\log\log\Delta})^2.
\end{align*}\\

\noindent{\bf Theorem 1.2\cite{ZD}}
If $G$ is a graph with the maximum degree $\Delta,$ then
\begin{align*}
\chi_{st}^{'}(G)\leq \chi_{st}^{'}(K_{\Delta+1})\cdot O(\frac{\log\Delta}{\log\log\Delta})^2,
\end{align*}
therefore,
\begin{align*}
\chi_{st}^{'}(G)\leq\Delta\cdot2^{O(1)\sqrt{\log\Delta}}.
\end{align*}

\noindent{\bf Theorem 1.3\cite{ZD}}

 (a) When $G$ is a subcubic graph, then $\chi^{\prime}_{st}(G)\leq7.$

 (b) If $G$ is a simple cubic graph, then $\chi^{\prime}_{st}(G)\geq4,$ and the equality holds if and only if it covers the graph of the 3-cubic.\\

\noindent{\bf Conjecture 1.4\cite{ZD}}
Let $G$ be a subcubic graph, then $\chi^{\prime}_{st}(G)\leq6.$\\

Bezegov$\mathrm{\acute{a}}$ et al.  have given some results and a conjecture on outerplanar graph in \cite{LB} as follows:\\

\noindent{\bf Theorem 1.5\cite{LB}}
If $G$ is an outerplanar graph, then\\

(1)  $\chi_{st}^{'} (G) \leq \lfloor 1.5\Delta \rfloor + 12.$

(2)If $\Delta(G)\leq3,$ then $\chi_{st}^{'} (G)  \leq 5 .$\\

\noindent{\bf Conjecture 1.6\cite{LB}}
If $G$ is an outerplanar graph with the maximum degree $\Delta$, then $\chi_{st}^{'} (G) \leq \lfloor 1.5\Delta \rfloor + 1.$\\

The following theorems 1.7-1.10 give some results on the star chromatic indexes of outerplanar graphs.\\

\noindent{\bf Theorem 1.7\cite{LKD}}
When $G$ is a maximal outerplanar graph with maximum degree $\Delta=4$, then $\chi_{st}^{'}(G)= 6.$\\

\noindent{\bf Theorem 1.8\cite{KD1}}
Let $G$ be a maximal outerplanar graph with $n$ vertices, then $6\leq\chi_{st}^{'} (G) \leq n-1.$\\

\noindent{\bf Theorem 1.9\cite{KD1}}
If $G$ is a maximal outerplanar graph with maximum degree $\Delta $ and $n(n\geq 8)$ vertices, then $\chi_{st}^{'}(F_{\Delta+1})\leq\chi_{st}^{'} (G) \leq n.$\\

Y. Q. Wang  et al. obtained the following result by edge partition method.\\

\noindent{\bf Theorem 1.10\cite{WF}}
If $G$ is an outerplanar graph with the maximum degree $\Delta,$  then $\chi_{st}^{'} (G) \leq \lfloor 1.5\Delta \rfloor + 5.$\\

\noindent{\bf Lemma 1.11\cite{GC}}
A graph $G$ is an outerplanar graph if and only if $G$ is $K_4$ and $K_{2,3}$ minor free.\\

In \cite{BB}, one has proved that almost all graphs have diameter 2, so it is important to discuss the star edge coloring of outerplanar graphs with  diameters 2 or 3.

For the sake of narrative, let $\zeta_n^d= \{G_{n} | G_{n}$ is an 2-connected outerplanar graph with diameter $d$  and $n$ vertices$\}$. In this paper, we give some upper bound based on maximum degree of graphs  $G_n(\in \zeta_n^2 \cup \zeta_n^3).$\\

\section{$\chi^{\prime}_{st}(G_n(\in \zeta_n^2))$}
\vspace{2mm}
First of all, we recall the star chromatic indexes of three basic graphs\cite{LKD}.\\

\noindent{\bf Lemma 2.1\cite{LKD}} For the path $P_n$ with $(n \geq5)$ vertices, we have $\chi_{st}^{'}(P_n)=3.$\\

\noindent{\bf Lemma 2.2\cite{LKD}}
Let $C_n$ be a cycle with $n$ vertices, then
\begin{align}
\chi^{\prime}_{st}(C_n)=
\begin{cases}
4,&  \;  n=5 , \\
3,&  \text{else}.
\end{cases}
\end{align}\\

\noindent{\bf Lemma 2.3\cite{LKD}}
  $F_{n+1}$ is a Fan graph with $n+1$ vertices, i.e. $F_{n+1}=K_1 \vee P_n,$  then
\begin{align}
\chi^{\prime}_{st}(F_{n+1})=
\begin{cases}
n+2,&  \;  n=4 , \\
n+1,&  \;  n=2,3,5,6 , \\
n,&    \; {n\geq7}.
\end{cases}
\end{align}\\

\noindent{\bf Definition 2.4}
When $G$ is a planar graph, we say that $G$ is a maximal planar graph if $G$ is a nonplanar graph when adding an edge for any two nonadjacent vertices.\\

\noindent{\bf Definition 2.5}
Let $G$ be an outerplanar graph.  When adding an edge for any two nonadjacent vertices, the resulting graph is a nonouterplanar graph, then $G$ is called a maximal outerplanar graph.\\


We know that any outerplanar graph can become a maximal outerplanar graph by adding edges from the definition of the maximal outerplanar graph, and the diameter of the resulting graph will not be increased. Thus, we only need to study the star chromatic index of maximal outerplanar graph with  diameter 1 and  diameter 2 in order to study the upper bound of the star chromatic index of 2-connected outerplanar graph with diameter 2. Since only $K_3$ has  diameter 1 and the connected subgraph with diameter 2 of $K_3$ is $P_3,$ which has no star edge coloring. Therefore, we wouldn't consider the maximal outerplanar with diameter 1.

In \cite{LR}, L. Beineke and R. Pippert gave a method to construct the maximal outerplanar graph shown as follows:

(a) $K_3$ is a maximal outerplanar graph;

(b) Let $G_1$ be a maximal outerplanar graph embedded in a plane with vertices lying on the exterior face $F_1$  and $G_2$ be the resulting graph by adding a new vertex which is adjacent to some two vertices of an edge on $F_1.$ Then $G_2$ is a maximal outerplanar graph;

(c) Any maximal outerplanar graph can be obtained by the processes (a)(b).\\

\noindent{\bf Theorem 2.6}
Let $G_n \in \zeta_n^2$, then
\begin{align}
\Delta(G_{n}) \le\chi^{\prime}_{st}(G_{n})\le
\begin{cases}
n+1,&  \;  n=5 , \\
n,&  \;  n=4,6,7 , \\
\Delta(G_{n}),&    \; {n\geq8}.
\end{cases}
\end{align}

\textbf{Proof:} By the above construction method, we know that maximal outerplanar graphs $G_n$ with diameter 2 have the following two cases: (1) $G_n=F_n$($n \geq 4$), (2) $G_n=G_6^{1}$.

\begin{figure}[htpb]
\centering
\begin{tikzpicture}[scale=0.5]
\tikzstyle{every node}=[font=\small,scale=0.5]
\fill (2,0) circle(3pt) ;
\node [below=2pt] at (2,0) {$v_1$};
\fill (5,0) circle(3pt) ;
\node [below=2pt] at (5,0) {$v_0$};
\fill (8,0) circle(3pt) ;
\node [below=2pt] at (8,0) {$v_4$};
\fill (5,5.1961) circle(3pt) ;
\node [above=2pt] at (5,5.1961) {$v_5$};
\fill (3.5,2.59805) circle(3pt) ;
\node [left=2pt] at (3.5,2.59805) {$v_2$};
\fill (6.5,2.59805) circle(3pt) ;
\node [right=2pt] at (6.5,2.59805) {$v_3$};
\draw (2,0)--node[below=0.5pt]{$1$}(5,0)--node[below=0.5pt]{$4$}(8,0)--node[right=0.5pt]{$6$}(6.5,2.59805)--node[right=0.5pt]{$4$}(5,5.1961)
--node[left=0.5pt]{$1$}(3.5,2.59805)--node[left=0.5pt]{$6$}(2,0);
\draw (5,0)--node[right=0.5pt]{$3$}(6.5,2.59805)--node[above=0.5pt]{$5$}(3.5,2.59805)
--node[left=0.5pt]{$2$}(5,0);
\fill (-2,0) circle(3pt) ;
\node [below=2pt] at (-2,0) {$v_4$};
\fill (-5,0) circle(3pt) ;
\node [below=2pt] at (-5,0) {$v_0$};
\fill (-8,0) circle(3pt) ;
\node [below=2pt] at (-8,0) {$v_1$};
\fill (-5,5.1961) circle(3pt) ;
\node [above=2pt] at (-5,5.1961) {$v_5$};
\fill (-3.5,2.59805) circle(3pt) ;
\node [right=2pt] at (-3.5,2.59805) {$v_3$};
\fill (-6.5,2.59805) circle(3pt) ;
\node [left=2pt] at (-6.5,2.59805) {$v_2$};
\draw (-2,0)--(-5,0)--(-8,0)--(-6.5,2.59805)--(-5,5.1961)--(-3.5,2.59805)--(-2,0);
\draw (-5,0)--(-6.5,2.59805)--(-3.5,2.59805)--(-5,0);
\end{tikzpicture}
\caption{$G_6^{1}$}\label{F1}
\end{figure}
For the case (1), the Lemma 2.3 shows that the theorem is correct. In the case (2), $F_5$ is a vertex-induced subgraph of $G_6^{1},$ then $\chi_{st}^{'}(G_6^{1})\geq \chi_{st}^{'}(F_5)=6.$  It is easy to check that the coloring of $G_6^{1}$ shown in Fig. \ref{F1} is a star 6-edge coloring.

According to the definition of the star edge coloring, we have $\chi_{st}^{'}(G_{n}) \geq \Delta$. For $n \geq 8$, we know that $\Delta(G_n)=n-1$. $\hfill \Box$


\section{$\chi^{\prime}_{st}(G_n(\in \zeta_n^3))$}

\vspace{2mm}
By the definition of maximal outerplanar graphs and Lemma 1.11,  we only need to study the star chromatic index of the maximal outerplanar graph with diameter 2 or 3 for studying the upper bound of the star chromatic index of $G_n.$ Let $\xi_n^{d}= \{H_{n} | H_{n}$ are the 2-connected maximal outerplanar graph with diameter $d$ and $n$ vertices $\}.$
  Let $n$ be an integer, for any $G_n,$ there is some $H_n\in \xi _n^{2} \cup \xi _n^{3}$ such that $G_n$ is a subgraph of $H_n$. The following two cases are discussed according to the diameter of $H_n$.

\subsection{$H_n\in \xi _n^{2}$}
\noindent{\bf Theorem 3.1}
Let $H_n\in \xi_n^2$ and  $G_n$ is an edge-include subgraph of $H_n.$ When $H_n=F_n$, then we have

$$\Delta(G_{n}) \leq \chi_{st}^{'}(G_{n})\leq \Delta(G_{n})+3.$$

If $H_n=G_6^{1}$, we obtained that
\begin{equation}
\chi_{st}^{'}(G_{n})=4.
\end{equation}
\textbf{Proof:} In the following, we prove the theorem by two cases based on $H_n$ with diameter 2.
\subsubsection{$H_n=G_6^{1}$}
We delete any two edges of $\{v_0v_2,v_0v_3,v_2v_3\}$ in $H_n$, then the diameter of resulting graph is 3. In the sense of  isomorphism, we only consider the graph ${G_6^{1}}^{'}$ shown in Fig. \ref{F2}, which is obtained by deleting the edges  $v_0v_2,v_0v_3$ in $H_n$.
\begin{figure}[htpb]
\centering
\begin{tikzpicture}[scale=0.5]
\tikzstyle{every node}=[font=\small,scale=0.5]
\fill (0,0) circle(3pt) ;
\node [left=2pt] at (0,0) {$v_1$};
\fill (4,0) circle(3pt) ;
\node [right=2pt] at (4,0) {$v_4$};
\fill (4,2) circle(3pt) ;
\node [right=2pt] at (4,2) {$v_3$};
\fill (0,2) circle(3pt) ;
\node [left=2pt] at (0,2) {$v_2$};
\fill (2,3.59805) circle(3pt) ;
\node [above=2pt] at (2,3.59805) {$v_5$};
\fill (2,-1.59805) circle(3pt) ;
\node [below=2pt] at (2,-1.59805) {$v_0$};
\draw (0,0)--(2,-1.59805)--(4,0)--(4,2)--(2,3.59805)--(0,2)--(0,0);
\draw (0,2)--(4,2);
\fill (7,0) circle(3pt) ;
\node [left=2pt] at (7,0) {$v_1$};
\fill (11,0) circle(3pt) ;
\node [right=2pt] at (11,0) {$v_4$};
\fill (11,2) circle(3pt) ;
\node [right=2pt] at (11,2) {$v_3$};
\fill (7,2) circle(3pt) ;
\node [left=2pt] at (7,2) {$v_2$};
\fill (9,3.59805) circle(3pt) ;
\node [above=2pt] at (9,3.59805) {$v_5$};
\fill (9,-1.59805) circle(3pt) ;
\node [below=2pt] at (9,-1.59805) {$v_0$};
\draw (7,0)--node[left=0.5pt]{$1$}(9,-1.59805)--node[right=0.5pt]{$2$}(11,0)--node[right=0.5pt]{$3$}(11,2)--node[right=0.5pt]{$2$}(9,3.59805)--node[left=0.5pt]{$1$}(7,2)--node[left=0.5pt]{$3$}(7,0);
\draw (7,2)--node[above=0.5pt]{$4$}(11,2);
\end{tikzpicture}
\caption{${G_6^{1}}^{'}$}\label{F2}
\end{figure}
Since the edges $v_1v_2, v_2v_5, v_5v_3, v_3v_4$ must be colored by at least three different colors, and there are no  bichromatic 4-length path in ${G_6^{1}}^{'}$, moreover the color of $v_2v_3$ is different from the colors of $v_1v_2, v_2v_3, v_5v_3, v_3v_4$, thus we obtain that $\chi_{st}^{'}({G_{6}^1}^{'})\geq 4$. It is easy to check that the coloring of Fig. \ref{F2} is a star edge coloring and only uses four colors, then $\chi_{st}^{'}({G_{6}^1}^{'})=4$.\\

\subsubsection{$H_n=F_n$} Let $H_n=F_n$, then $n \geq 6$(since diam($C_5$)=2).
We discuss according to the value of $n$ as follows.\\

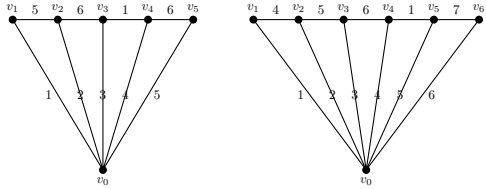
\begin{figure}[htpb]
\centering
\begin{tikzpicture}[scale=0.5]
\tikzstyle{every node}=[font=\small,scale=0.5]
\fill (0,0) circle(3pt) ;
\node [below=2pt] at (0,0) {$v_0$};
\fill (0,4) circle(3pt) ;
\node [above=2pt] at (0,4) {$v_3$};
\fill (-1.2,4) circle(3pt) ;
\node [above=2pt] at (-1.2,4) {$v_2$};
\fill (-2.4,4) circle(3pt) ;
\node [above=2pt] at (-2.4,4) {$v_1$};
\fill (1.2,4) circle(3pt) ;
\node [above=2pt] at (1.2,4) {$v_4$};
\fill (2.4,4) circle(3pt) ;
\node [above=2pt] at (2.4,4) {$v_5$};
\draw (0,0)--node[left=0.5pt]{$1$}(-2.4,4)--node[above=0.5pt]{$5$}(-1.2,4)--node{$2$}(0,0)--node{$3$}(0,4)--node[above=0.5pt]{$1$}(1.2,4)--node{$4$}(0,0)--node[right=0.5pt]{$5$}(2.4,4)--node[above=0.5pt]{$6$}(1.2,4);
\draw (0,4)--node[above=0.5pt]{$6$}(-1.2,4);
\fill (7,0) circle(3pt) ;
\node [below=2pt] at (7,0) {$v_0$};
\fill (6.4,4) circle(3pt) ;
\node [above=2pt] at (6.4,4) {$v_3$};
\fill (5.2,4) circle(3pt) ;
\node [above=2pt] at (5.2,4) {$v_2$};
\fill (4,4) circle(3pt) ;
\node [above=2pt] at (4,4) {$v_1$};
\fill (7.6,4) circle(3pt) ;
\node [above=2pt] at (7.6,4) {$v_4$};
\fill (8.8,4) circle(3pt) ;
\node [above=2pt] at (8.8,4) {$v_5$};
\fill (10,4) circle(3pt) ;
\node [above=2pt] at (10,4) {$v_6$};
\draw (7,0)--node[left=0.5pt]{$1$}(4,4)--node[above=0.5pt]{$4$}(5.2,4)--node{$2$}(7,0)--node{$3$}(6.4,4)--node[above=0.5pt]{$6$}(7.6,4)--node{$4$}(7,0)--node{$5$}(8.8,4)--node[above=0.5pt]{$1$}(7.6,4);
\draw (6.4,4)--node[above=0.5pt]{$5$}(5.2,4);
\draw (7,0)--node[right=0.5pt]{$6$}(10,4)--node[above=0.5pt]{$7$}(8.8,4);
\end{tikzpicture}
\caption{$F_6$ and $F_7$}\label{F3}
\end{figure}

{\bf Case 1.} If $n=6,$  the star edge coloring of $H_6$ is shown in Fig. \ref{F3}. When we delete the edge(s) which is(are) not in the exterior face of $H_6$, the resulting graph has diameter 3. Let the set of deleting edge(s) to be $D$. We can obtain that $\chi_{st}^{'}(G_{6})=\chi_{st}^{'}(H_{6})-\vert D\vert \leq 6-\vert D\vert=\Delta(G_{6})+1.$\\

{\bf Case 2.} When $n=7$, the star edge coloring of $H_n$ is shown in Fig. \ref{F3}. Similarly,
when we delete the edge(s) $D$  in the interior face of $H_7$, the diameter of the resulting graph becomes 3. Therefore, we obtain that $\chi_{st}^{'}(G_{7})\leq \chi_{st}^{'}(C_{7})+\vert D\vert=3+\vert D\vert=\Delta(G_{7})+1.$

{\bf Case 3.} If $n\geq 8$, to ensure $diam(G_{n})=3$, by enumeration, we can obtain that the deleted edges $D$ ($\vert D\vert\leq3$). Therefore, we can get the formula as follows, $\chi_{st}^{'}(G_{n})\leq\chi_{st}^{'}(H_{n})=\Delta(H_{n})\leq\Delta(G_{n})+3.$
 \\

%

Finally, when $G_n$ is an edge induced subgraph of $H_n=F_n(\in \xi_n^2)$, we obtain that
 $\Delta(G_{n}) \leq \chi_{st}^{'}(G_{n})\leq \Delta(G_{n})+3.$ $\hfill \Box$

\subsection{$H_n\in \xi_n^3$}
We firstly discuss the property of $H_n$ in the following.\\

 \noindent {\bf Theorem~3.2}\ \ If $H_n\in \xi_n^3,$  then $H_n$ has $G_6^1$ or $G_6^2$(shown in the Fig. 4) as a vertex induced subgraph. \\

 \textbf{Proof:} By the construction method of maximal outerplanar graph in \cite{LR},  maximal outerplanar graph of order 5 is $F_5$ under the meaning of isomorphism, thus $F_5$ must be a vertex induced subgraph of $H_n(n\geq 5).$  Assume that the maximum degree vertex of $F_5$ is $v_0$ and other vertices are $v_1,v_2,$ $v_3,v_4.$  We add a vertex $v_5$ to $F_5$ in order to obtain a big maximal outerplanar graph: if $v_5$ is adjacent to $v_i$, $v_{i+1}(i\in \{1,2,3\}),$ the theorem is proven; if $v_5$ is adjacent to $v_0$,$v_1$ or $v_0$,$v_4,$ we obtain a Fan $F_6,$ then, we continue to add vertex, there is some vertex $v_j$ which is adjacent to $v_i$,$v_{i+1}(i \in \{1,2,3,\cdots ,j-1\}),$  otherwise, the diameter of the resulting graph is 2. $\hfill \Box$
\begin{figure}[htpb]
\begin{minipage}[H]{0.5\linewidth}
\centering
\begin{tikzpicture}[scale=0.5]
\tikzstyle{every node}=[font=\small,scale=0.5]
\fill (0,0.5) circle(3pt) ;
\node [below=2pt] at (0,0.5) {$v_0$};
\fill (0,4) circle(3pt) ;
\node [above=2pt] at (0,4) {$v_2$};
\fill (-1.5,4) circle(3pt) ;
\node [above=2pt] at (-1.5,4) {$v_1$};
\fill (1.5,4) circle(3pt) ;
\node [above=2pt] at (1.5,4) {$v_3$};
\fill (3,4) circle(3pt) ;
\node [above=2pt] at (3,4) {$v_4$};
\fill (2.2,5.8667) circle(3pt) ;
\node [above=2pt] at (2.2,5.8667) {$v_5$};
\draw (0,0.5)--(-1.5,4)--(0,4)--(0,0.5)--(1.5,4)--(2.2,5.8667)--(3,4)--(0,0.5);
\draw (0,4)--(1.5,4)--(3,4);
\end{tikzpicture}
\caption{$G_6^2$}\label{F4}
\end{minipage}
\begin{minipage}[H]{0.5\linewidth}
\centering
\begin{tikzpicture}[scale=0.5]
\tikzstyle{every node}=[font=\small,scale=0.5]
\fill (0,0) node (v0) {} circle(3pt) ;
\fill (-3,0) node (v1) {} circle(3pt) ;
\fill (-1.5,2.599) node (v2) {} circle(3pt) ;
\fill (1.5,2.599) node (v3) {} circle(3pt) ;
\fill (3,0) node (v4) {} circle(3pt) ;
\fill (0,5.196) node (v5) {} circle(3pt) ;
\draw (v0)--(v1)--(v2)--(v5)--(v3)--(v4)--(v0);
\draw (v0)--(v2)--(v3)--(v0);
\node [above=3pt] at (0,0) {$v_0$};
\node [right=2pt] at (-3,0.1) {$v_1$};
\node [right=2pt] at (3,0.1) {$v_4$};
\node [right=2pt] at (-1.5,2.699) {$v_2$};
\node [left=2pt] at (1.5,2.699) {$v_3$};
\node [above=2pt] at (0,5.196) {$v_5$};
\fill (2.61,-0.67) node (v6) {} circle(3pt) ;
\fill (2.2,-1.2) node (v7) {} circle(3pt) ;
\fill (1.4,-1.68) node (v8) {} circle(3pt) ;
\draw (v0)--(v8)--(v7)--(v0)--(v6)--(v7);
\draw (v4)--(v6);
\fill (3.5,0.8) node (v9) {} circle(3pt) ;
\fill (3.6,1.5) node (v10) {} circle(3pt) ;
\draw (v4)--(v9)--(v3)--(v10)--(v9);
\fill (3.4,2.899) node (v11) {} circle(3pt) ;
\fill (2.9,3.699) node (v12) {} circle(3pt) ;
\fill (2.15,4.499) node (v13) {} circle(3pt) ;
\fill (1.2,5) node (v14) {} circle(3pt) ;
\draw (v5)--(v14)--(v3)--(v13)--(v14);
\draw (v3)--(v11)--(v12)--(v3);
\draw (v12)--(v13);
\fill (-0.9,5.096) node (v15) {} circle(3pt) ;
\fill (-1.795,4.896) node (v16) {} circle(3pt) ;
\fill (-2.5,4.496) node (v17) {} circle(3pt) ;
\draw (v5)--(v15)--(v2)--(v16)--(v15);
\draw (v2)--(v17)--(v16);
\fill (-3.4,2.799) node (v18) {} circle(3pt) ;
\fill (-3.5,2.1) node (v19) {} circle(3pt) ;
\fill (-3.55,1.4) node (v20) {} circle(3pt) ;
\fill (-3.3,0.7) node (v21) {} circle(3pt) ;
\draw (v2)--(v18)--(v19)--(v2)--(v20)--(v19);
\draw (v1)--(v21)--(v2);
\draw (v21)--(v20);
\fill (-2.7,-0.8) node (v22) {} circle(3pt) ;
\fill (-1.6,-1.6) node (v23) {} circle(3pt) ;
\fill (-0.4,-1.78) node (v24) {} circle(3pt) ;
\draw (v1)--(v22)--(v0)--(v23)--(v22);
\draw (v0)--(v24)--(v23);
\end{tikzpicture}
\caption{$H_n(G_6^1 \subseteq H_n)$}
\end{minipage}
\end{figure}

\noindent {\bf Theorem~3.3}\ \ Let $G$ be a connected graph.  When $G$ has two adjacent vertices $v_0$,$v_1\in G$ and the distance between $v_0$ and any other vertex of $G$ is less than $diam(G)$, then adding some vertices  $v_2,v_3,\cdots v_n$ and edges $v_1v_2, v_2v_3, \cdots ,v_{n-1}v_n, v_0v_2, v_0v_3, \cdots ,v_0v_n$, we obtain a new graph $G^{'}$  with $diam(G^{'})$ $=diam(G).$ \\

  By theorem 3.3, we know that $v_0v_1$,$v_1v_2$,$v_2v_5$,$v_5v_3$,$v_3v_4$,$v_4v_0$ can be extended to be Fans as shown in Fig. 5. Assume that $\Delta (H_n)=\Delta,$ then $\Delta =max\{ d(v_0),$
$d(v_2),d(v_3)\}$. Especially,  the graph of $d(v_0)=d(v_2)=d(v_3)= \Delta $ is denoted by $H_{n_\Delta},$ we give a lower bound of $\chi_{st}^{'}(H_{n_\Delta})$ in the following theorem.\\

\noindent {\bf Theorem~3.4}\ \ $\chi_{st}^{'}(H_{n_\Delta}) \geq \Delta(H_{n_\Delta}) +2$.\\

\textbf{Proof:}  We consider the subgraph $G_{n_\Delta}$ (shown as the Fig. 6)  of $H_{n_\Delta}.$   We firstly prove $\chi_{st}^{'}(G_{n_\Delta}) \geq \Delta +1.$ It is easy to know that $\chi_{st}^{'}(G_{n_\Delta}) \geq \Delta$ by the definition of star edge coloring.  Assume that
$\chi_{st}^{'}(G_{n_\Delta}) = \Delta,$ since $d(v_0)=\Delta$, there are some vertex $v_0^{'}\in N(v_0)$ such that $c(v_0v_0^{'})=c(v_2v_3)$. Similarly, there are some vertex $v_3^{'}\in N(v_3)$ such that $c(v_3v_3^{'})=c(v_0v_2)$, thus the path $v_0^{'}-v_0-v_2-v_3-v_3^{'}$ is 2-edge colored, which contradicts to the definition of star edge coloring. Therefore, we obtain that $\chi_{st}^{'}(H_{n_\Delta})\geq \chi_{st}^{'}(G_{n_\Delta}) \geq \Delta +1.$\\
\indent Now, we assume that $\chi_{st}^{'}(G_{n_\Delta}) = \Delta +1$. Without loss of generality, let $c(v_0v_1)=1,$  $c(v_0v_0^{(1)})=2,$  $c(v_0v_0^{(2)})=3,$ $\cdots,$  $c(v_0v_0^{(\Delta -4)})=\Delta -3,$ $c(v_0v_4)=\Delta -2,$ $c(v_0v_3)=\Delta -1,$  $c(v_0v_2)=\Delta,$  $c(v_2v_3)=\Delta +1$
(If there are some vertex $v_{x_1} \in \{ v_1, v_0^{(1)},v_0^{(2)},$ $\cdots v_0^{(\Delta-4)},v_4 \},$  such that $c(v_0v_{x_1})=c(v_2v_3),$  then for any vertex $v_{y} \in \{ v_1, v_2^{(1)},v_2^{(2)},\cdots v_2^{(\Delta-4)},v_5 \}$, we have $c(v_2v_y)\neq c(v_0v_3).$ The coloring scheme is isomorphic to the form: $c(v_0v_1)=1,$
$c(v_0v_0^{(1)})=2,$ $c(v_0v_0^{(2)})=3,$  $\cdots,$  $c(v_0v_0^{(\Delta -4)})=\Delta -3,$ $c(v_0v_4)=\Delta -2,$  $c(v_0v_3)=\Delta -1,$  $c(v_0v_2)=\Delta,$  $c(v_2v_3)=\Delta +1$.\\

Let $\{c(v_3v_4), c(v_3v_5),c(v_3v_3^{(1)}),c(v_3v_3^{(2)}),\cdots, c(v_3v_3^{(\Delta-4)})\}=Q.$  Firstly, we prove that $\Delta-2\notin Q:$
since $c(v_3v_4)\neq \Delta+1~\mbox{or}~\Delta-2,$ there is some vertex $v_{x_2}\in\{v_2,v_1,v_0^{(1)},v_0^{(2)},\cdots,v_0^{(\Delta-4)}\}$ such that
$c(v_0v_{x_2})=c(v_3v_4),$ which forces that there is a 4-path colored by $\Delta-2,c(v_3v_4)$ in $G_{n_\Delta}$ when $\Delta-2\in Q.$   The coloring is not proper.  \\
\indent Since $\{\Delta- 2,\Delta-1,\Delta+1\}\notin Q$ and $\left| Q \right|=\Delta- 2$, $Q=\{ 1,2,\cdots,\Delta -3,\Delta \}$. Similarly, we have $\{1,\Delta,\Delta+1\}\notin \{c(v_2v_1),c(v_2v_5),c(v_2v_2^{(1)}),c(v_2v_2^{(2)}), \cdots, c(v_2v_2^{(\Delta-4)}) \}$. Since $\Delta\in Q$, $\Delta-1\notin \{c(v_2v_1),c(v_2v_5)$,$c(v_2v_2^{(1)}),c(v_2v_2^{(2)}), \cdots, c(v_2v_2^{(\Delta-4)}) \}$, a contradiction.
$\hfill \Box$\\
\\
{\bf Theorem~3.5}
If $H_n\in \xi_n^3,$ and $G_n \in \zeta_n^{3},$ then $\Delta(G_{n}) \leq \chi_{st}^{'}(G_{n}) \leq \Delta(G_{n})+6;$    when $G_n \in \xi_n^{3}$,  we have $\Delta(G_{n}) \leq \chi_{st}^{'}(G_{n}) \leq \Delta(G_{n})+4.$\\

\textbf{Proof:} We prove the theorem through the following two cases 3.5.1 and 3.5.2.\\


\textbf{Case 3.5.1}\ \ $H_n$ has $G_6^1$ as an induced subgraph. \\


\indent We get the Fig. 7 (denoted by $H^{'}_{n_\Delta}$) by adding edges $v_1v_0^{(1)},$ $v_0^{(i)}v_0^{(i+1)},$ $v_2^{(1)}v_1$,
$v_2^{(j)}v_2^{(j+1)},$ $v_2^{(\Delta-4)}v_5$, $v_5v_3^{(1)},$ $v_3^{(k)}v_3^{(k+1)},$ $v_3^{(\Delta-4)}v_4$,
 $v_0^{(\Delta-4)}v_4$ in $G_{n_\Delta}$ where $1\leq i,j,k\leq \Delta-3.$  In the following, we discuss the upper bound of $\chi_{st}^{'}(H^{'}_{n_\Delta})$ in terms of $\Delta$.
\begin{figure}[htpb]
\begin{minipage}[H]{0.5\linewidth}
\centering
\begin{tikzpicture}[scale=0.5]
\tikzstyle{every node}=[font=\small,scale=0.5]
\fill (0,0) node (v0) {} circle(3pt) ;
\fill (-3,0) node (v1) {} circle(3pt) ;
\fill (-1.5,2.599) node (v2) {} circle(3pt) ;
\fill (1.5,2.599) node (v3) {} circle(3pt) ;
\fill (3,0) node (v4) {} circle(3pt) ;
\fill (0,5.196) node (v5) {} circle(3pt) ;
\draw (v0)--(v1)--(v2)--(v5)--(v3)--(v4)--(v0);
\draw (v0)--(v2)--(v3)--(v0);
\node [above=3pt] at (0,0) {$v_0$};
\node [right=2pt] at (-3,0.1) {$v_1$};
\node [right=2pt] at (3,0.1) {$v_4$};
\node [right=2pt] at (-1.5,2.699) {$v_2$};
\node [left=2pt] at (1.5,2.699) {$v_3$};
\node [above=2pt] at (0,5.196) {$v_5$};
\fill (2.4,-1.2) node (v6) {} circle(3pt) ;
\fill (1.7,-1.68) node (v7) {} circle(3pt) ;
\node [right=3pt] at (2.4,-1.1) {$v_0^{(\Delta -4)}$};
\node [right=3pt] at (1.7,-1.68) {$v_0^{(\Delta -5)}$};
\draw (v0)--(v6);
\draw (v0)--(v7);
\fill (3.5,1.1) node (v8) {} circle(3pt) ;
\fill (3.6,2.2) node (v9) {} circle(3pt) ;
\node [right=3pt] at (3.5,1.1) {$v_3^{(\Delta -4)}$};
\node [right=3pt] at (3.6,2.2) {$v_3^{(\Delta -5)}$};
\draw (v3)--(v8);
\draw (v3)--(v9);
\fill (2.1,3.55) node (v30) {} circle(1pt) ;
\fill (2.3,3.25) node (v30) {} circle(1pt) ;
\fill (2.45,2.89) node (v30) {} circle(1pt) ;
\fill (2.15,4.499) node (v10) {} circle(3pt) ;
\fill (1.2,5) node (v11) {} circle(3pt) ;
\node [above=2pt] at (2.35,4.499) {$v_3^{(2)}$};
\node [above=2pt] at (1.4,5) {$v_3^{(1)}$};
\draw (v3)--(v11);
\draw (v3)--(v10);
\fill (-0.9,5.096) node (v12) {} circle(3pt) ;
\fill (-1.895,4.796) node (v13) {} circle(3pt) ;
\node [above=2pt] at  (-0.9,5.096) {$v_2^{(\Delta -4)}$};
\node [above=2pt] at (-1.895,4.796) {$v_2^{(\Delta -5)}$};
\draw (v12)--(v2);
\draw (v13)--(v2);
\fill (-2.1,3.55) node (v30) {} circle(1pt) ;
\fill (-2.3,3.25) node (v30) {} circle(1pt) ;
\fill (-2.45,2.89) node (v30) {} circle(1pt) ;
\fill (-3.5,2.2) node (v14) {} circle(3pt) ;
\fill (-3.55,1.1) node (v15) {} circle(3pt) ;
\node [left=2pt] at (-3.5,2.2) {$v_2^{(2)}$};
\node [left=2pt] at  (-3.55,1.1) {$v_2^{(1)}$};
\draw (v2)--(v14);
\draw (v2)--(v15);
\fill (-2.35,-1.2) node (v16) {} circle(3pt) ;
\fill (-1.2,-1.78) node (v17) {} circle(3pt) ;
\node [left=2pt] at (-2.35,-1.2) {$v_0^{(1)}$};
\node [left=2pt] at  (-1.2,-1.78) {$v_0^{(2)}$};
\draw (v0)--(v16);
\draw (v0)--(v17);
\fill (-0.4,-1.15) node (v30) {} circle(1pt) ;
\fill (0,-1.25) node (v30) {} circle(1pt) ;
\fill (0.45,-1.19) node (v30) {} circle(1pt) ;
\end{tikzpicture}
\caption{$G_{n_\Delta}$} \label{t1}
\end{minipage}
\begin{minipage}[H]{0.5\linewidth}
\centering
\begin{tikzpicture}[scale=0.5]
\tikzstyle{every node}=[font=\small,scale=0.5]
\fill (0,0) node (v0) {} circle(3pt) ;
\fill (-3,0) node (v1) {} circle(3pt) ;
\fill (-1.5,2.599) node (v2) {} circle(3pt) ;
\fill (1.5,2.599) node (v3) {} circle(3pt) ;
\fill (3,0) node (v4) {} circle(3pt) ;
\fill (0,5.196) node (v5) {} circle(3pt) ;
\draw (v0)--(v1)--(v2)--(v5)--(v3)--(v4)--(v0);
\draw (v0)--(v2)--(v3)--(v0);
\node [above=3pt] at (0,0) {$v_0$};
\node [right=2pt] at (-3,0.1) {$v_1$};
\node [right=2pt] at (3,0.1) {$v_4$};
\node [right=2pt] at (-1.5,2.699) {$v_2$};
\node [left=2pt] at (1.5,2.699) {$v_3$};
\node [above=2pt] at (0,5.196) {$v_5$};
\fill (2.61,-0.67) node (v6) {} circle(3pt) ;
\fill (2.2,-1.2) node (v7) {} circle(3pt) ;
\fill (1.4,-1.68) node (v8) {} circle(3pt) ;
\draw (v0)--(v8)--(v7)--(v0)--(v6)--(v7);
\draw (v4)--(v6);
\node [below=2pt] at (1.49,-1.68)  {$v_0^{(i+1)}$};
\node [below=2pt] at (-0.4,-1.78)  {$v_0^{(i)}$};
\fill (3.5,0.8) node (v9) {} circle(3pt) ;
\fill (3.6,1.5) node (v10) {} circle(3pt) ;
\draw (v4)--(v9)--(v3)--(v10)--(v9);
\node [right=2pt] at (3.6,1.5)  {$v_3^{(k+1)}$};
\node [right=2pt] at (3.4,2.899)  {$v_3^{(k)}$};
\fill (3.4,2.899) node (v11) {} circle(3pt) ;
\fill (2.9,3.699) node (v12) {} circle(3pt) ;
\fill (2.15,4.499) node (v13) {} circle(3pt) ;
\fill (1.2,5) node (v14) {} circle(3pt) ;
\draw (v5)--(v14)--(v3)--(v13)--(v14);
\draw (v3)--(v11)--(v12)--(v3);
\draw (v12)--(v13);
\node [left=2pt] at (-3.4,2.799)  {$v_2^{(j)}$};
\node [left=2pt] at (-2.3,4.696)  {$v_2^{(j+1)}$};
\fill (-0.9,5.096) node (v15) {} circle(3pt) ;
\fill (-1.795,4.896) node (v16) {} circle(3pt) ;
\fill (-2.5,4.496) node (v17) {} circle(3pt) ;
\draw (v5)--(v15)--(v2)--(v16)--(v15);
\draw (v2)--(v17)--(v16);
\fill (-3.4,2.799) node (v18) {} circle(3pt) ;
\fill (-3.5,2.1) node (v19) {} circle(3pt) ;
\fill (-3.55,1.4) node (v20) {} circle(3pt) ;
\fill (-3.3,0.7) node (v21) {} circle(3pt) ;
\draw (v2)--(v18)--(v19)--(v2)--(v20)--(v19);
\draw (v1)--(v21)--(v2);
\draw (v21)--(v20);
\fill (-2.7,-0.8) node (v22) {} circle(3pt) ;
\fill (-1.6,-1.6) node (v23) {} circle(3pt) ;
\fill (-0.4,-1.78) node (v24) {} circle(3pt) ;
\draw (v1)--(v22)--(v0)--(v23)--(v22);
\draw (v0)--(v24)--(v23);
\draw[dashed]  (3.6,1.5)--(3.4,2.899);
\draw[dashed]  (-0.4,-1.78)--(1.49,-1.68) ;
\draw[dashed] (-3.4,2.799)--(-2.3,4.696) ;
\end{tikzpicture}
\caption{$H^{'}_{n_\Delta}$}\label{t2}
\end{minipage}
\end{figure}\\
  (1) When $\Delta=5,$ the Fig. 8(a) gives a star edge coloring, then $\Delta+2=7 \le \chi_{st}^{'}(H^{'}_{n_\Delta}) \le 9= \Delta+4$.\\
  (2) If $\Delta=6,$ a star edge coloring is given in the Fig. 8(b), we obtain that $\Delta+2=8 \le \chi_{st}^{'}(H^{'}_{n_\Delta}) \le 9=\Delta+3 $.\\
  (3) For the case $\Delta=7,$ it is easy to check that the coloring of the Fig. 8(c) is a star edge coloring, thus $\Delta+2=9 \le \chi_{st}^{'}(H^{'}_{n_\Delta}) \le 11= \Delta+4 $.\\
  (4) $\Delta=8$, we consider the star edge coloring of the Fig. 8(d), thus $\Delta+2=10 \le \chi_{st}^{'}(H^{'}_{n_\Delta}) \le 11=\Delta+3 $.\\
  (5) When $\Delta \geq 9,$  we consider the following coloring function $c:$ $c(v_0v_1)=1$, $c(v_0v_4)=\Delta-2$, $c(v_2v_1)=2$, $c(v_2v_5)=\Delta+2$, $c(v_3v_5)=1$, $c(v_3v_4)=\Delta+2$, $c(v_0v_2)=\Delta$, $c(v_0v_3)=\Delta-1$, $c(v_2v_3)=\Delta+1$, $c(v_0v_0^{(l)})=l+1$, $c(v_2v_2^{(l)})=l+2$, $c(v_3v_3^{(l)})=l+1$, $c(v_1v_0^{(1)})=\Delta+3$, $c(v_0^{(m-1)}v_0^{(m)})=m+3$, $c(v_0^{(\Delta -5)}v_0^{(\Delta -4)})=\Delta +1$, $c(v_0^{(\Delta -4)}v_4)=\Delta+3$, $c(v_1v_2^{(1)})=5$, $c(v_2^{(n-1)}v_2^{(n)})=n+4$,  $c(v_2^{(\Delta -5)}v_2^{(\Delta -4)})=2$, $c(v_2^{(\Delta -4)}v_5)=3$, $c(v_5v_3^{(1)})=4$, $c(v_3^{(p-1)}v_3^{(p)})=p+3$,  $c(v_3^{(\Delta-6)}v_3^{(\Delta-5)})=\Delta +2$, $c(v_3^{(\Delta-5)}v_3^{(\Delta-4)})=\Delta$, $c(v_3^{(\Delta-4)}v_4)=2$, where $l,m,p \in Z_{+}$ and $m \in [2,\Delta-5]$, $n \in [2,\Delta-5]$, $p \in [2,\Delta-6]$. As an example, the Fig.8(e) gives a star edge coloring of $H^{'}_{n_9}.$

 Next, we prove that the coloring $c$ is a star edge coloring: the paths $v_1-v_2^{(1)}-v_2^{(2)}-\cdots -v_2^{(\Delta-5)}-v_2^{(\Delta-4)}-v_5,$ $v_5-v_3^{(1)}-v_3^{(2)}-\cdots -v_3^{(\Delta-5)}-v_3^{(\Delta-4)}-v_4,$ $v_1-v_0^{(1)}-v_0^{(2)}-\cdots -v_0^{(\Delta-5)}-v_0^{(\Delta-4)}-v_4$ are denoted by $P_2,$ $P_3$, $P_0,$ respectively. By the Theorem 1.4 of $^{[5]},$ we know that the Fans $P_0 \vee v_0$, $P_2 \vee v_2$, $P_3 \vee v_3$ has no bichromatic 4-paths(cycles).

 Now, assume $P$ is a  bichromatic 3-path(cycle)of $H^{'}_{n_\Delta}.$  If the path(cycle)$P$ contains $v_0v_2$, $v_0v_3$ or $v_2v_3$ as an edge, then all the bichromatic 3-path(cycle) $P$ of $H^{'}_{n_\Delta}$ are $v_2^{(i)}-v_2-v_3-v_3^{(i+1)}$,$v_1-v_2-v_3-v_3^{(1)}$, $v_5-v_2-v_3-v_4$, $v_2^{(i)}-v_2-v_0-v_0^{(i+1)}$, $v_1-v_2-v_0-v_0^{(1)}$, $v_2^{(5)}-v_2-v_0-v_4$, $v_0^{(j)}-v_0-v_3-v_3^{(j)}$, $v_1-v_0-v_3-v_5$,where $i,j \in Z_{+}$ and $i \in [1,\Delta-5]$, $j \in [1,\Delta-4]$. It is easy to check that any one of the bichromatic 3-path(cycle) can't be extended to become bichromatic 4-paths(cycles), thus $H^{'}_{n_\Delta}$ has no bichromatic 4-paths(cycles) containing $v_0v_2$, $v_0v_3$ or $v_2v_3$ as an edge. \\
\indent If the path(cycle) $P$ contains $v_1v_0^{(1)}$, $v_1v_2^{(1)}$, $v_4v_0^{(\Delta-4)}$, $v_4v_3^{(\Delta-4)}$, $v_5v_2^{(\Delta-4)}$ or $v_5v_3^{(1)}$ as an edge,  then all the bichromatic 3-path(cycle) $P$ of $H^{'}_{n_\Delta}$ are(the path $P$ containing $v_0v_2$, $v_0v_3$ or $v_2v_3$ as an edge can be discussed similarly) $v_2^{(1)}-v_1-v_2-v_2^{(3)}$, $v_1-v_2^{(1)}-v_2-v_2^{(3)}$, $v_2^{(1)}-v_1-v_0-v_0^{(4)}$, $v_2^{(\Delta-4)}-v_5-v_2-v_2^{(1)}$, $v_5-v_2^{(\Delta-4)}-v_2-v_2^{(1)}$, $v_2^{(\Delta-4)}-v_5-v_3-v_3^{(2)}$, $v_3^{(1)}-v_5-v_3-v_3^{(3)}$, $v_5-v_3^{(1)}-v_3-v_3^{(3)}$, $v_3^{(1)}-v_5-v_2-v_2^{(2)}$, $v_3^{(\Delta-4)}-v_4-v_3-v_3^{(1)}$, $v_4-v_3^{(\Delta-4)}-v_3-v_3^{(1)}$, $v_3^{(\Delta-4)}-v_4-v_0-v_0^{(1)}$.  It is easy to check that any one of the bichromatic 3-path(cycle) can't be extended to become bichromatic 4-paths(cycles), thus the graph $H^{'}_{n_\Delta}$ has no bichromatic 4-paths(cycles) containing $v_1v_0^{(1)}$, $v_1v_2^{(1)}$, $v_4v_0^{(\Delta-4)}$, $v_4v_3^{(\Delta-4)}$, $v_5v_2^{(\Delta-4)}$ or $v_5v_3^{(1)}$ as an edge. \\
\indent In the following, we consider the paths $P$ containing $v_0v_1$, $v_1v_2$, $v_2v_5$, $v_5v_3$, $v_3v_4$ or $v_4v_0$ as an edge. All the bichromatic 3-path(cycle) $P$ of $H^{'}_{n_\Delta}$ are(the path $P$ containing $v_0v_2$, $v_0v_3$, $v_2v_3$, $v_1v_0^{(1)}$, $v_1v_2^{(1)}$, $v_4v_0^{(5)}$, $v_4v_3^{(5)}$, $v_5v_2^{(5)}$ or $v_5v_3^{(1)}$ as an edge can be discussed similarly)  $v_1-v_2-v_2^{(\Delta-4)}-v_2^{(\Delta-5)}$, $v_1-v_2-v_2^{(\Delta-5)}-v_2^{(\Delta-4)}$, $v_2-v_1-v_0-v_0^{(1)}$, $v_2-v_5-v_3-v_4$, $v_4-v_3-v_3^{(\Delta-6)}-v_3^{(\Delta-5)}$, $v_4-v_3-v_3^{(\Delta-5)}-v_3^{(\Delta-6)}$, $v_4-v_0-v_0^{(\Delta-6)}-v_0^{(\Delta-5)}$, $v_4-v_0-v_0^{(\Delta-5)}-v_0^{(\Delta-6)}$.  One can find that any one of the bichromatic 3-path(cycle) can't be extended to become bichromatic 4-paths(cycles), thus the graph $H^{'}_{n_\Delta}$ has no bichromatic 4-paths(cycles) containing $v_0v_1$, $v_1v_2$, $v_2v_5$, $v_5v_3$, $v_3v_4$ or $v_4v_0$ as an edge.

 Let $H^ {''}_{n_\Delta}$ be the subgraph by deleting $v_0v_2$, $v_0v_3$, $v_2v_3$, $v_1v_0^{(1)}$, $v_1v_2^{(1)}$, $v_4v_0^{(\Delta-4)}$, $v_4v_3^{(\Delta-4)}$, $v_5v_2^{(\Delta-4)}$, $v_5v_3^{(1)}$, $v_0v_1$, $v_1v_2$, $v_2v_5$, $v_5v_3$, $v_3v_4$ and $v_4v_0$ in $H^{'}_{n_\Delta}.$  The graph $H^{''}_{n_\Delta}$ is not connected and has three connected components, which are vertex-induced subgraphs of Fans $P_0\vee v_0$, $P_2\vee v_2$  and $P_3\vee v_3,$ thus $H^{''}_{n_\Delta}$  doesn't contian bichromatic 4-paths(cycles), therefore $H^{'}_{n_\Delta}$ has no bichromatic 4-paths(cycles).

 So overall, we obtain that $\Delta+2 \leq \chi_{st}^{'}(H_{n_\Delta}) \leq \Delta+3$ when $\Delta\geq 9.$
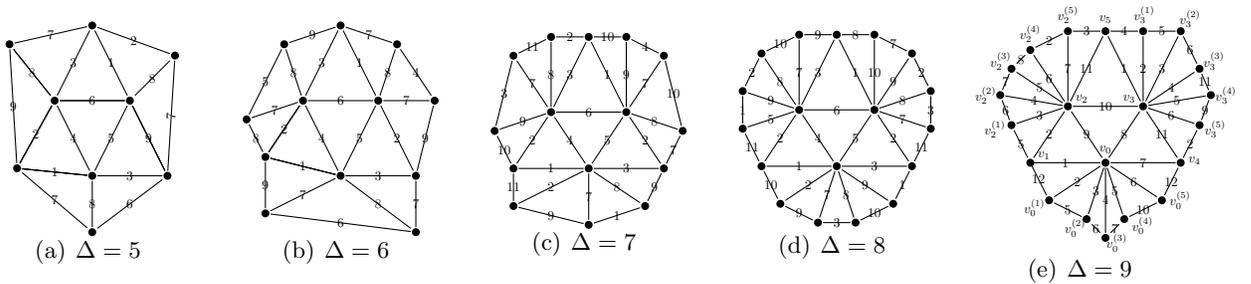
\begin{figure}[htpb]
\centering
\subfigure[$\Delta=5$]{
\begin{minipage}{3.04cm}
\centering
\begin{tikzpicture}[scale=0.5]
\tikzstyle{every node}=[font=\small,scale=0.5]
\fill (0,-2) circle(3pt) ;
\fill (-2,-1.8) circle(3pt) ;
\fill (-2.2,1.5) circle(3pt) ;
\fill (2,-2) circle(3pt) ;
\fill (2.2,1.2) circle(3pt) ;
\fill (1,0) circle(3pt) ;
\fill (-1,0) circle(3pt) ;
\fill (0,2) circle(3pt) ;
\fill (0,-3.5) circle(3pt) ;
\node (v3) at (0,-2) {};
\node (v2) at (0,-3.5) {};
\node (v4) at (-2,-1.8) {};
\node (v5) at (-2.2,1.5) {};
\node (v8) at (2,-2) {};
\node (v12) at (0,2) {};
\node (v10) at (2.2,1.2) {};
\node (v7) at (1,0) {};
\draw (v4) -- node{$1$}(v3) --(v4) --node{$9$} (v5);
\draw (v3) --node{$8$}(v2) --node{$7$}(v4);
\draw (v2) --node{$6$}(v8);
\draw (v10) --node[pos=0.5,sloped]{$7$}(v8);
\node (v6) at (-1,0) {};
\draw (v5) -- (v6) --node{$2$} (v4) --(v6) -- node{$4$}(v3)
-- (v3) -- node{$5$}(v7) -- (v6);
\draw (v3) -- node{$3$}(v8) ;
\draw (v7) -- node{$9$}(v8)  -- node{$$}(v7) -- node{$8$}(v10);
\draw (v5)--node{$7$}(v12);
\draw (v5) -- node{$8$}(v6) -- node{$3$}(v12) -- node{$1$}(v7) ;
\draw (v12) -- node{$2$}(v10);
\draw (v6) --node{$6$}(v7);

\end{tikzpicture}
\end{minipage}}
\subfigure[$\Delta=6$]{
\begin{minipage}{3.04cm}
\centering
\begin{tikzpicture}[scale=0.5]
\tikzstyle{every node}=[font=\small,scale=0.5]
\fill (0,-2) circle(3pt) ;
\fill (-2,-3) circle(3pt) ;
\fill (-2,-1.5) circle(3pt) ;
\fill (-2.5,-0.5) circle(3pt) ;
\fill (2,-2) circle(3pt) ;
\fill (2.5,0) circle(3pt) ;
\fill (1.5,1.5) circle(3pt) ;
\fill (1,0) circle(3pt) ;
\fill (-1,0) circle(3pt) ;
\fill (0,2) circle(3pt) ;
\fill (-1.5,1.5) circle(3pt) ;
\fill (2,-3.5) circle(3pt) ;
\node (v3) at (0,-2) {};
\node (v1) at (-2,-3) {};
\node (v2) at (2,-3.5) {};
\node (v4) at (-2,-1.5) {};
\node (v5) at (-2.5,-0.5) {};
\node (v8) at (2,-2) {};
\node (v9) at (2.5,0) {};
\node (v10) at (1.5,1.5) {};
\node (v7) at (1,0) {};
\draw (v1) --node{$6$} (v2) -- node{$8$}(v3) -- node{$7$}(v1) --node{$9$} (v4)
-- node{$1$}(v3) --(v4) --node{$8$} (v5);
\node (v6) at (-1,0) {};
\draw (v5) --node{$7$} (v6) --node{$2$} (v4) --(v6) -- node{$4$}(v3)
-- (v3) -- node{$5$}(v7) -- node{$6$}(v6);
\draw (v3) -- node{$3$}(v8)  -- node{$7$}(v2) ;
\draw (v7) -- node{$2$}(v8)  -- node{$9$}(v9) -- node{$7$}(v7) -- node{$8$}(v10)
-- node{$4$}(v9);
\node (v12) at (0,2) {};
\node (v11) at (-1.5,1.5) {};
\draw (v5) -- node{$5$}(v11) -- node{$8$}(v6) -- node{$3$}(v12) -- node{$1$}(v7) ;
\draw (v11) -- node{$9$}(v12) -- node{$7$}(v10);
\end{tikzpicture}
\end{minipage}}
\subfigure[$\Delta=7$]{
\begin{minipage}{3.04cm}
\centering
\begin{tikzpicture}[scale=0.5]
\tikzstyle{every node}=[font=\small,scale=0.5]
\tikzstyle{every node}=[font=\small,scale=0.5]
\fill (0,2) circle(3pt) ;
\node (v12) at (0,2) {};
\node (v11) at (-2,1.5) {};
\fill (-2,1.5) node (v21) {} circle(3pt) ;
\node at (0,-3) {};
\fill (0,-3) node (v16) {} circle(3pt) ;
\node at (-1,2) {};
\fill (-1,2) node (v23) {} circle(3pt) ;
\node (v3) at (0,-1.5) {};
\fill (0,-1.5) node (v15) {} circle(3pt) ;
\node (v1) at (-2,-2.5) {};
\fill (-2,-2.5) node (v13) {} circle(3pt) ;
\node (v2) at (1.5,-2.5) {};
\fill (1.5,-2.5) node (v17) {} circle(3pt) ;
\node (v4) at (-2,-1.5) {};
\fill (-2,-1.5) node (v14) {} circle(3pt) ;
\node (v5) at (-2.5,-0.5) {};
\fill (-2.5,-0.5) node (v19) {} circle(3pt) ;
\node (v8) at (2,-1.5) {};
\fill (2,-1.5) node (v18) {} circle(3pt) ;
\node (v9) at (2.5,-0.5) {};
\fill (2.5,-0.5) node (v26) {} circle(3pt) ;
\node (v10) at (2,1.5) {};
\fill (2,1.5) node (v25) {} circle(3pt) ;
\node (v7) at (1,0) {};
\fill (1,0) node (v22) {} circle(3pt) ;
\node (v6) at (-1,0) {};
\fill (-1,0) node (v20) {} circle(3pt) ;
\node at (1,2) {};
\fill (1,2) node (v24) {} circle(3pt) ;
\draw (v13) --node{$11$} (v14) --node{$1$} (v15) --node{$2$} (v13)
--node{$9$} (v16) --node{$1$} (v17) -- node{$8$}(v15) -- node{$3$}(v18) -- node{$9$}(v17);
\draw (v15) --node{$7$} (v16);
\draw (v19) -- node{$10$} (v14) -- node{$2$} (v20) --node{$9$}  (v19)
-- node{$3$} (v21) --node{$7$}  (v20) -- node{$4$} (v15) -- node{$5$} (v22)
-- node{$6$} (v20) --node{$8$}  (v23) -- node{$11$} (v21);
\draw (v23) -- node{$2$} (v12) -- node{$10$} (v24) --node{$4$}  (v25)
--node{$10$}(v26) --   node{$7$} (v18) --node{$2$} (v22) -- node{$7$} (v25);
\draw (v20) -- node{$3$}(v12) -- node{$1$}(v22) -- node{$9$}(v24);
\draw (v22) --node{$8$} (v26) ;

\end{tikzpicture}
\end{minipage}}
\subfigure[$\Delta=8$]{
\begin{minipage}{3.04cm}
\centering
\begin{tikzpicture}[scale=0.5]
\tikzstyle{every node}=[font=\small,scale=0.5]
\fill (0,2) node (v31) {}circle(3pt) ;
\node (v12) at (0,2) {};
\node (v11) at (-2,1.5) {};
\fill (-2,1.5) node (v21) {} circle(3pt) ;
\node at (-0.5,-3) {};
\fill (-0.5,-3) node (v16) {} circle(3pt) ;
\node at (-1,2) {};
\fill (-1,2) node (v23) {} circle(3pt) ;
\node (v3) at (0,-1.5) {};
\fill (0,-1.5) node (v15) {} circle(3pt) ;
\node (v2) at (1.5,-2.5) {};
\fill (1.5,-2.5) node (v17) {} circle(3pt) ;
\node (v4) at (-2,-1.5) {};
\fill (-2,-1.5) node (v14) {} circle(3pt) ;
\node (v5) at (-2.5,-0.5) {};
\fill (-2.5,-0.5) node (v19) {} circle(3pt) ;
\node (v8) at (2,-1.5) {};
\fill (2,-1.5) node (v18) {} circle(3pt) ;
\node (v9) at (2.5,-0.5) {};
\fill (2.5,-0.5) node (v32) {} circle(3pt) ;
\node (v10) at (2,1.5) {};
\fill (2,1.5) node (v25) {} circle(3pt) ;
\node (v7) at (1,0) {};
\fill (1,0) node (v22) {} circle(3pt) ;
\node (v6) at (-1,0) {};
\fill (-1,0) node (v20) {} circle(3pt) ;
\node at (1,2) {};
\fill (1,2) node (v24) {} circle(3pt) ;
\node at (0.5,-3) {};
\fill (0.5,-3) node (v26) {} circle(3pt) ;
\node at (-2.5,0.5) {};
\fill (-2.5,0.5) node (v29) {} circle(3pt) ;
\node at (2.5,0.5) {};
\fill (2.5,0.5)node (v30) {}  circle(3pt) ;
\node at (-1.5,-2.5) {};
\fill (-1.5,-2.5) node (v1) {} circle(3pt) ;
\draw (v1) -- node{$10$}(v14) --node{$1$} (v15) --node{$2$} (v1) --node{$9$} (v16) --
node{$7$}(v15) -- node{$8$}(v26) --node{$3$}(v16);
\draw (v26) -- node{$10$}(v17) -- node{$9$}(v15) -- node{$3$}(v18) --node{$1$}(v17);
\draw (v15) -- node{$5$}(v22) -- node{$2$}(v18);
\draw (v18) -- node{$11$}(v32) -- node{$7$}(v22) --node{$8$} (v30) -- node{$3$}(v32);
\draw (v30) -- node{$2$}(v25) -- node{$9$}(v22) -- node{$6$}(v20) -- node{$5$}(v19)
-- node{$11$}(v14) -- node{$2$}(v20) --node{$4$} (v15);
\draw (v25) -- node{$7$}(v24) -- node{$10$}(v22) -- node{$1$}(v12) -- node{$3$}(v20)
--  node{$9$}(v29) --  node{$1$}(v19);
\draw (v24) --node{$8$}(v12) -- node{$9$}(v23) -- node{$7$}(v20) -- node{$8$}(v21)
-- node{$2$}(v29);
\draw (v21) -- node{$10$}(v23);
\end{tikzpicture}
\end{minipage}}
\subfigure[$\Delta=9$]{
\begin{minipage}{3.04cm}
\centering
\begin{tikzpicture}[scale=0.5]
\tikzstyle{every node}=[font=\small,scale=0.5]
\fill (0,2) node (v31) {}circle(3pt) ;
\node [above=2pt] at (0,2) {$v_5$};
\fill (-2,1.5) node (v21) {} circle(3pt) ;
\node [above=2pt] at (-2,1.5) {$v_2^{(4)}$};
\fill (-0.5,-3) node (v16) {} circle(3pt) ;
\node [below=2pt] at (-0.8,-2.8) {$v_0^{(2)}$};
\fill (-1,2) node (v23) {} circle(3pt) ;
\node [above=2pt] at (-1,2) {$v_2^{(5)}$};
\fill (0,-1.5) node (v15) {} circle(3pt) ;
\node [above=2pt] at (0,-1.45) {$v_0$};
\fill (1.5,-2.5) node (v17) {} circle(3pt);
\node [right=2pt] at (1.5,-2.5) {$v_0^{(5)}$};
\fill (-2,-1.5) node (v14) {} circle(3pt) ;
\node [right=2pt] at (-2,-1.3) {$v_1$};
\fill (-2.5,-0.5) node (v19) {} circle(3pt) ;
\node [left=2pt] at (-2.45,-0.6) {$v_2^{(1)}$};
\fill (2,-1.5) node (v18) {} circle(3pt) ;
\node [right=2pt] at (2,-1.45) {$v_4$};
\fill (2.5,-0.5) node (v32) {} circle(3pt) ;
\node [right=2pt] at (2.43,-0.6) {$v_3^{(5)}$};
\fill (2,2) node (v25) {} circle(3pt) ;
\node [above=2pt] at (2.25,1.9) {$v_3^{(2)}$};
\fill (1,0) node (v22) {} circle(3pt) ;
\node [left=2pt] at (1,0.2) {$v_3$};
\fill (-1,0) node (v20) {} circle(3pt) ;
\node [right=2pt] at (-1,0.2) {$v_2$};
\fill (1,2) node (v24) {} circle(3pt) ;
\node [above=2pt] at (1.03,2) {$v_3^{(1)}$};
\fill (0,-3.5) node (v26) {} circle(3pt) ;
\node [below=2pt] at (0.3,-3.2) {$v_0^{(3)}$};
\fill (-2.5,1) node (v29) {} circle(3pt) ;
\node [left=2pt] at (-2.3,1.2) {$v_2^{(3)}$};
\fill (2.5,1)node (v30) {}  circle(3pt) ;
\node [right=2pt] at (2.4,1.2) {$v_3^{(3)}$};
\fill (-1.5,-2.5) node (v1) {} circle(3pt) ;
\node [left=2pt] at (-1.4,-2.7) {$v_0^{(1)}$};
\fill (-2.8,0.3) node (v36) {} circle(3pt) ;
\node [left=2pt] at (-2.7,0.3) {$v_2^{(2)}$};
\fill (2.8,0.3) node (v33) {} circle(3pt) ;
\node [right=2pt] at (2.72,0.3) {$v_3^{(4)}$};
\fill (0.5,-3) node (v34) {} circle(3pt) ;
\node [right=2pt] at (0.5,-3.17) {$v_0^{(4)}$};
\draw (v1) -- node{$12$}(-2,-1.5) -- node (v2) {$1$}(v15) -- node{$2$}(v1)
-- node{$5$}(v16) -- node{$3$}(v15) -- node{$4$}(v26) -- node{$6$}(v16);
\draw (v26) -- node{$7$}(v34) --node{$5$} (v15) --node{$6$} (v17) --node{$10$} (v34);
\draw (v17) -- node{$12$}(v18) -- node{$7$}(v15) --node{$8$} (v22) -- node{$11$}(v18)
-- node{$2$}(v32) -- node{$6$}(v22) -- node{$5$}(v33) -- node{$9$}(v32);
\draw (v33) -- node{$11$}(v30) -- node{$4$}(v22) -- node{$10$}(v20) -- node{$9$}(v15);
\draw (v30) -- node{$6$}(v25) -- node{$3$}(v22) -- node{$1$}(v31) --node{$4$} (v24)
-- node{$5$}(v25);
\draw (v24) -- node{$2$}(v22);
\draw (v31) -- node{$3$}(v23) --node{$7$} (v20) --node{$6$} (v21) -- node{$8$}(v29)
-- node{$5$}(v20) -- node{$11$}(v31);
\draw (v21) -- node{$2$}(v23);
\draw (v29) --node{$7$} (v36) -- node{$4$}(v20) -- node{$3$}(v19) -- node{$6$}(v36);
\draw (v20) -- node{$2$}(v14)--node{$5$}(v19);
\end{tikzpicture}
\end{minipage}}
\caption{ Star edge coloring of $H^{'}_{n_\Delta}$ ($\Delta =5,6,\cdots ,9$)} \label{t1}
\end{figure}
  When $\Delta=5~\mbox{or}~ 7,$  we proved that $\chi_{st}^{'}(H_{n_\Delta}) \leq \Delta+4$. In the end, we have $\Delta \leq \chi_{st}^{'}(H_{n}) \leq \chi_{st}^{'}(H_{n_\Delta}) \leq \Delta+3$ when $\Delta\geq 8~\mbox{or} ~\Delta=6.$ $\hfill \Box$\\

\textbf{Case 3.5.2}\ \ If $H_n$ does not contain $G_6^1,$ then $H_n$ must contain $G_6^2$ as an edge induced subgraph. \\

\indent By the construction method of maximal outerplanar graph in \cite{LR}, it is easy to check that there is no vertex adjacent to $v_0,$ $v_4$ or $v_2,$  $v_3,$ otherwise $H_n$ must contain  an edge induced subgraph $G_6^1.$  Since $H_n$ has diameter three,
 $H_n$  has two cases shown as Fig. 9.\\

\begin{figure}[htpb]
\begin{minipage}[H]{0.5\linewidth}
\centering
\begin{tikzpicture}[scale=0.45]
\tikzstyle{every node}=[font=\small,scale=0.45]
\fill (-1.5,1.5) node (v1) {} circle(3pt) ;
\node [below=2pt] at (-1.5,1.5) {$v_3$};
\fill (-6,1) node (v2) {} circle(3pt) ;
\node [left=1pt] at (-6,1) {$v_0^{(\Delta -4)}$};
\fill (-7,-1) node (v3) {} circle(3pt) ;
\node [left=1pt] at (-7,-1) {$v_0^{(3)}$};
\fill (-7.25,-2) node (v4) {} circle(3pt) ;
\node [left=1pt] at (-7.25,-2) {$v_0^{(2)}$};
\fill (-7.5,-3) node (v5) {} circle(3pt) ;
\node [left=1pt] at (-7.5,-3) {$v_0^{(1)}$};
\fill (-2,-3.5) node (v6) {} circle(3pt) ;
\node [below=2pt] at (-2,-3.5) {$v_0$};
\fill (-3.5,1.5) node (v8) {} circle(3pt) ;
\node [below=2pt] at (-3.5,1.5) {$v_2$};
\fill (1.5,1.5) node (v9) {} circle(3pt) ;
\node [below=2pt] at (1.5,1.5) {$v_4$};
\fill (-6.5,0.2) node (v10) {} circle(3pt) ;
\node [left=1pt] at (-6.5,0.2) {$v_0^{(\Delta -5)}$};
\fill (-5,1.5) node (v11) {} circle(3pt) ;
\node [below=2pt] at (-5,1.5) {$v_1$};
\fill (-3.1,3.1) node (v12) {} circle(3pt) ;
\node [left=1pt] at (-3.1,3.1) {$v_3^{(3)}$};
\fill (-3.3,2.6) node (v13) {} circle(3pt) ;
\node [left=1pt] at (-3.3,2.6) {$v_3^{(2)}$};
\fill (-3.5,2.1) node (v14) {} circle(3pt) ;
\node [left=1pt] at (-3.5,2.1) {$v_3^{(1)}$};
\draw (v1)--(v12)--(v13)--(v1)--(v14)--(v13);
\fill (-2.3,4.1) node (v15) {} circle(3pt);
\node [above=1pt] at (-2.3,4.1) {$v_3^{(\Delta -5)}$};
\fill (-1.5,4.4) node (v16) {} circle(3pt);
\fill (-0.5,4.5) node (v17) {} circle(3pt);
\node [above=2pt] at (-0.5,4.5) {$v_5$};
\draw (v1)--(v15)--(v16)--(v1)--(v17)--(v16);
\fill (-6.5,-1) node (v56) {} circle(1pt);
\fill (-6.4,-0.8) node (v57) {} circle(1pt);
\fill (-6.3,-0.6) node (v58) {} circle(1pt);
\fill (-6.2,-0.4) node (v59) {} circle(1pt);
\fill (1.8,3.6) node (v57) {} circle(1pt);
\fill (2.05,3.45) node (v58) {} circle(1pt);
\fill (2.3,3.3) node (v59) {} circle(1pt);
\fill (-2.4,3.6) node (v57) {} circle(1pt);
\fill (-2.6,3.4) node (v58) {} circle(1pt);
\fill (-2.8,3.2) node (v59) {} circle(1pt);
\fill (3.5,2.1) node (v18) {} circle(3pt);
\node [right=1pt] at (3.5,2.1) {$v_4^{(\Delta -3)}$};
\fill (3.1,2.7) node (v19) {} circle(3pt);
\node [right=1pt] at (3.1,2.7) {$v_4^{(\Delta -4)}$};
\fill (2.7,3.3) node (v20) {} circle(3pt);
\node [right=1pt] at (2.7,3.3) {$v_4^{(\Delta -5)}$};
\draw (v9)--(v18)--(v19)--(v9)--(v20)--(v19);
\fill (1.6,4.2) node (v21) {} circle(3pt);
\node [right=1pt] at (1.6,4.2) {$v_4^{(3)}$};
\fill (1,4.5) node (v22) {} circle(3pt);
\node [above=1pt] at (1,4.5) {$v_4^{(2)}$};
\fill (0.2,4.6) node (v23) {} circle(3pt);
\node [above=1pt] at (0.2,4.6) {$v_4^{(1)}$};
\draw (v9)--(v17)--(v23)--(v9)--(v22)--(v23);
\draw (v22)--(v21)--(v9);
\draw (v6) -- (v5) -- (v4) -- (v6) -- (v3) -- (v4);
\draw (v6)--(v2)--(v10)--(v6);
\draw (v2) -- (v11)--(v6)--(v8)--(v11);
\draw (v6)--(v9)--(v1)--(v6);
\draw (v1)--(v8);
\end{tikzpicture}
\end{minipage}
\begin{minipage}[H]{0.5\linewidth}
\centering
\begin{tikzpicture}[scale=0.5]
\tikzstyle{every node}=[font=\small,scale=0.5]
\fill (0,0) node (v0) {} circle(3pt) ;
\node [below=2pt] at (0,0) {$v_0$};
\fill (-2.5,4) node (v2) {} circle(3pt) ;
\node [below=2pt] at (-2.5,4) {$v_2$};
\fill (-4,4) node (v1) {} circle(3pt) ;
\node [below=2pt] at (-4,4) {$v_1$};
\fill (2.5,4) node (v3) {} circle(3pt) ;
\node [below=2pt] at (2.5,4) {$v_3$};
\fill (4,4) node (v4) {} circle(3pt) ;
\node [below=2pt] at (4,4) {$v_4$};
\draw (v0)--(v1)--(v2)--(v0)--(v4)--(v3)--(v0);
\draw (v2)--(v3);
\fill (-3.25,7) node (v6) {} circle(3pt) ;
\node [above=2pt] at (-3.25,7) {$v_6$};
\fill (3.25,7) node (v5) {} circle(3pt) ;
\node [above=2pt] at (3.35,7) {$v_5$};
\draw (v1)--(v6)--(v2);
\draw (v4)--(v5)--(v3);
\fill (-2.7,7.2) node (v21) {} circle(3pt) ;
\node [above=2pt] at (-2.7,7.2) {$v_2^{(1)}$};
\fill (-2.1,7.15) node (v22) {} circle(3pt) ;
\node [above=2pt] at (-2.1,7.15) {$v_2^{(2)}$};
\fill (-1.6,7) node (v23) {} circle(3pt) ;
\node [above=2pt] at (-1.4,7) {$v_2^{(3)}$};
\draw (v6)--(v21)--(v2)--(v22)--(v21);
\draw (v2)--(v23)--(v22);
\fill (-1.1,6.5) node (v24) {} circle(3pt) ;
\node [right=1pt] at (-1.1,6.7) {$v_2^{(\Delta-6)}$};
\fill (-0.8,5.9) node (v25) {} circle(3pt) ;
\node [right=1pt] at (-0.8,6.1) {$v_2^{(\Delta-5)}$};
\fill (-0.5,5.3) node (v26) {} circle(3pt) ;
\node [right=1pt] at (-0.5,5.5) {$v_2^{(\Delta-4)}$};
\draw (v24)--(v2)--(v25)--(v24);
\draw (v2)--(v26)--(v25);
\fill (2.7,6.9) node (v30) {} circle(3pt) ;
\node [above=1pt] at (2.7,6.9) {$v_3^{(\Delta-4)}$};
\fill (2.2,6.55) node (v31) {} circle(3pt) ;
\fill (1.8,6.2) node (v32) {} circle(3pt) ;
\node [above=1pt] at (1.7,6.2) {$v_3^{(\Delta-6)}$};
\draw (v5)--(v30)--(v3)--(v31)--(v30);
\draw (v3)--(v32)--(v31);
\fill (1.1,5.4) node (v33) {} circle(3pt) ;
\fill (0.7,4.95) node (v34) {} circle(3pt) ;
\fill (0.3,4.5) node (v35) {} circle(3pt) ;
\node [left=1pt] at (0.3,4.6) {$v_3^{(1)}$};
\draw (v33)--(v3)--(v34)--(v33);
\draw (v3)--(v35)--(v34);
\fill (1.5,5.3) node (v40) {} circle(0.8pt) ;
\fill (1.65,5.5) node (v40) {} circle(0.8pt) ;
\fill (1.8,5.7) node (v40) {} circle(0.8pt) ;
\fill (-1.4,6.3) node (v40) {} circle(0.8pt) ;
\fill (-1.5,6.4) node (v40) {} circle(0.8pt) ;
\fill (-1.6,6.5) node (v40) {} circle(0.8pt) ;
\end{tikzpicture}
\end{minipage}
\caption{$H_n(G_6^2 \subseteq H_N)$}\label{t2}
\end{figure}
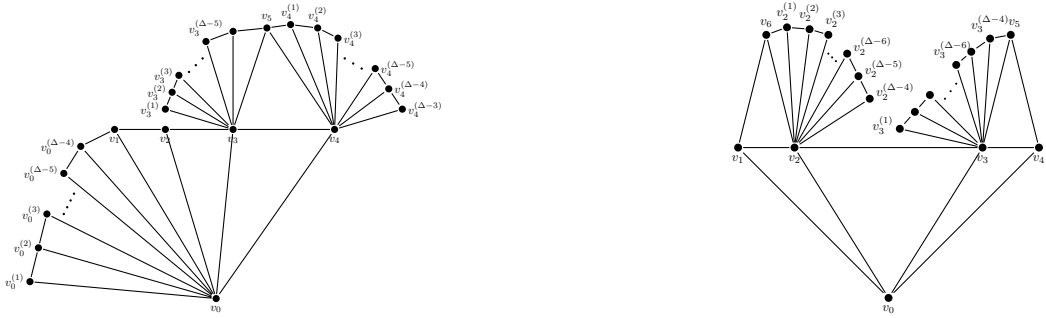

\noindent {\bf Case~1:}\ \ As shown in the left of Fig. 9, we discuss the case $d(v_0)=d(v_3)=d(v_4)=\Delta$ (denoted by $H_{n_\Delta}$). We consider the upper bound of
$\chi_ {st}^{'}(H_{n_\Delta})$ in terms of $\Delta.$ \\
 (1) If $\Delta =4,$ it is easy to check that the coloring of Fig. 10(a) is a star edge coloring with $\Delta +2$ colors, then $\chi_{st}^{'}(H_{n_\Delta}) \leq 6=\Delta +2.$  Since $F_5$ is an edge induced subgraph and $\chi_{st}^{'}(F_5)=\Delta +2$, thus $\chi_{st}^{'}(H_{n_\Delta})=\Delta +2$; \\
 (2)  The case of $\Delta =5,$  the Fig. 10(b) give a star edge coloring with $\Delta +3$ colors, therefore $\chi_{st}^{'}(H_{n_\Delta}) \leq 8=\Delta +3;$\\
 (3) For $\Delta =6$,  we consider the Fig. 10(c) which has a star edge coloring with $\Delta +3$ colors, thus $\chi_{st}^{'}(H_{n_\Delta}) \leq9 =\Delta +3$;\\
(4) When $\Delta \geq 7,$  we consider the following coloring function $c:$ $c(v_0v_1)=\Delta -3$, $c(v_0v_2)=\Delta -2$, $c(v_0v_3)=\Delta +1$, $c(v_0v_4)=\Delta +2$, $c(v_3v_5)=\Delta -1$, $c(v_4v_5)=\Delta +4$, $c(v_0v_0^{(l)})=l$, $c(v_3v_3^{(m)})=m+2$, $c(v_4v_4^{(p)})=p+2$, $c(v_1v_2)=\Delta+3 $, $c(v_2v_3)=\Delta$, $c(v_3v_4)=\Delta +3$, $c(v_0^{(l)}v_0^{(l+1)})=l+3$, $c(v_3^{(m)}v_3^{(m+1)})=m$, $c(v_4^{(p)}v_4^{(p+1)})=p$, $c(v_3^{(\Delta -4)}v_5)=\Delta -4$, $c(v_5v_4^{(1)})=\Delta,$ where $l,~m,~p \in Z_{+}$ and~ $l,~ m \in [l,\Delta-4]$, $p \in [l,\Delta-3]$.\\
\indent We prove the coloring function $c$ is a star edge coloring: the paths $v_0^{(1)}-v_0^{(2)}-\cdots -v_0^{(\Delta -4)}-v_1-v_2-v_3-v_4,$ $v_3^{(1)}-v_3^{(2)}-\cdots -v_3^{(\Delta -4)}-v_5,$ $v_5-v_4^ {(1)}-v_4^{(2)}-\cdots -v_4^{(\Delta -4)}-v_4^{(\Delta -3)}$ are denoted by $P_0$, $P_3,$ $P_4$ respectively. By the theorem 1.4 of $^{[5]},$  we know that the Fans $P_0 \vee v_0$, $P_3 \vee v_3$, $P_4 \vee v_4$ have no bichromatic paths or cycles of length four. Assume $P$ is a path(cycle) of $H_{n_\Delta}$ containing $v_0v_3$,$v_3v_4$,$v_2v_3$  or $v_5v_4^{(1)}$ and is a bichromatic path(cycle) with length three. All  bichromatic 3-path(cycle)of $H_ {n_\Delta}$ are $v_1-v_2-v_3-v_4,$ $v_2-v_3-v_5-v_4^{(1)}, $  $v_0-v_2-v_3-v_3^{(3)},$ $v_3^{(i)}-v_3-v_4-v_4^{(i)}\ (1\leqslant i\leqslant \Delta-4),$ $v_5-v_3-v_4-v_4^{(\Delta-3)},$ $v_4-v_4^{(1)}-v_5-v_3^{(3)}.$  By the coloring function $c,$ we know that all the paths(cycles) can't be become  bichromatic 4-paths(cycles) by extending,   thus $H_{n_\Delta}$ does not contain  bichromatic 4-paths(cycles) including $v_0v_3$,$v_3v_4$,$v_2v_3$ or $v_5v_4^{(1)}.$ Since  $c(v_0v_4)=\Delta+2$, $c(v_4v_5)=\Delta+4$, and $\Delta +2$, $\Delta +4$ appear only once in $H_{n_\Delta},$  deleting the edges $v_3v_4$, $v_2v_3$, $v_0v_3$, $v_0v_4$, $v_4v_5$ and $v_5v_4^{(1)}$, we obtain a subgraph denoted by $H^{'}_{n_\Delta}$. The subgraph $H^{'}_{n_\Delta}$ is not connected and has three connected components, which are vertex induced Fans $P_0\vee v_0$, $P_3 \vee v_3$ and $P_4\vee v_4,$  thus $H^{'}_{n_\Delta}$ has no bichromatic 4-paths and cycles. Therefore, $H_{n_\Delta}$ has no bichromatic 4-paths and cycles.  We obtain that $\chi_{st}^{'}(H_{n_\Delta}) \leq  \Delta+4$. As an example, the Fig. 10(d) gives the coloring function of $H_{n_7}.$


\begin{figure}[htpb]
\centering
\subfigure[$\Delta =4$]{
\begin{minipage}{3.8cm}
\centering
\begin{tikzpicture}[scale=0.52]
\tikzstyle{every node}=[font=\small,scale=0.52]
\fill (0,0) node (v0) {} circle(3pt) ;
\fill (0.6,3) node (v3) {} circle(3pt) ;
\fill (-0.6,3) node (v2) {} circle(3pt) ;
\fill (-1.9,3) node (v1) {} circle(3pt) ;
\fill (1.9,3) node (v4) {} circle(3pt) ;
\draw (v0)--node{$1$}(v1)--node{$4$}(v2)--node{$2$}(v0)--node{$3$}(v3)--node{$5$}(v2);
\draw (v0)--node{$4$}(v4)--node{$6$}(v3);
\fill (1.35,4.5) node (v5) {} circle(3pt) ;
\fill (2.3,4.8) node (v6) {} circle(3pt) ;
\draw (v3)--node{$1$}(v5)--node{$5$}(v4)--node{$3$}(v6)--node{$2$}(v5);
\end{tikzpicture}
\end{minipage}}
\subfigure[$\Delta =5$]{
\begin{minipage}{3.8cm}
\centering
\begin{tikzpicture}[scale=0.52]
\tikzstyle{every node}=[font=\small,scale=0.52]
\fill (0,0) node (v0) {} circle(3pt) ;
\fill (0.6,3) node (v3) {} circle(3pt) ;
\fill (-0.6,3) node (v2) {} circle(3pt) ;
\fill (-1.9,3) node (v1) {} circle(3pt) ;
\fill (1.9,3) node (v4) {} circle(3pt) ;
\draw (v0)--node{$2$}(v1)--node{$8$}(v2)--node{$3$}(v0)--node{$4$}(v3)--node{$1$}(v2);
\draw (v0)--node{$5$}(v4)--node{$8$}(v3);
\fill (1.35,4.5) node (v5) {} circle(3pt) ;
\fill (2.3,4.75) node (v6) {} circle(3pt) ;
\draw (v3)--node{$6$}(v5)--node{$3$}(v4)--node{$7$}(v6)--node{$1$}(v5);
\fill (-2.35,2.2) node (v7) {} circle(3pt) ;
\draw (v0)--node{$1$}(v7)--node{$5$}(v1);
\fill (0.45,4.4) node (v8) {} circle(3pt) ;
\fill (3,4.4) node (v9) {} circle(3pt) ;
\draw (v3)--node{$7$}(v8)--node{$2$}(v5);
\draw (v4)--node{$6$}(v9)--node{$2$}(v6);
\end{tikzpicture}
\end{minipage}}
\subfigure[$\Delta =6$]{
\begin{minipage}{3.8cm}
\centering
\begin{tikzpicture}[scale=0.52]
\tikzstyle{every node}=[font=\small,scale=0.52]
\fill (0,0) node (v0) {} circle(3pt) ;
\fill (0.6,3) node (v3) {} circle(3pt) ;
\fill (-0.6,3) node (v2) {} circle(3pt) ;
\fill (-1.9,3) node (v1) {} circle(3pt) ;
\fill (1.9,3) node (v4) {} circle(3pt) ;
\draw (v0)--node{$3$}(v1)--node{$8$}(v2)--node{$4$}(v0)--node{$5$}(v3)--node{$7$}(v2);
\draw (v0)--node{$6$}(v4)--node{$1$}(v3);
\fill (1.35,4.5) node (v5) {} circle(3pt) ;
\fill (2.2,4.75) node (v6) {} circle(3pt) ;
\draw (v3)--node{$3$}(v5)--node{$8$}(v4)--node{$2$}(v6)--node{$4$}(v5);
\fill (-2.5,2.4) node (v7) {} circle(3pt) ;
\draw (v0)--node{$2$}(v7)--node{$9$}(v1);
\fill (0.6,4.4) node (v8) {} circle(3pt) ;
\draw (v3)--node{$2$}(v8)--node{$9$}(v5);
\fill (2.9,4.4) node (v9) {} circle(3pt) ;
\draw (v4)--node{$4$}(v9)--node{$5$}(v6);
\fill (-3,1.6) node (v10) {} circle(3pt) ;
\draw (v0)--node{$1$}(v10)--node{$4$}(v7);
\fill (0,3.9) node (v11) {} circle(3pt) ;
\draw (v3)--node{$8$}(v11)--node{$4$}(v8);
\fill (3.3,3.7) node (v12) {} circle(3pt) ;
\draw (v4)--node{$7$}(v12)--node{$3$}(v9);
\end{tikzpicture}
\end{minipage}}
\subfigure[$\Delta =7$]{
\begin{minipage}{3.8cm}
\centering
\begin{tikzpicture}[scale=0.52]
\tikzstyle{every node}=[font=\small,scale=0.52]
\fill (0,0) node (v0) {} circle(3pt) ;
\node [below=2pt] at (0,0) {$v_0$};
\fill (0.6,3) node (v3) {} circle(3pt) ;
\node [below=2pt] at (0.8,3) {$v_3$};
\fill (-0.6,3) node (v2) {} circle(3pt) ;
\node [below=2pt] at (-0.3,3.01) {$v_2$};
\fill (-1.9,3) node (v1) {} circle(3pt) ;
\node [below=2pt] at (-1.6,3.03) {$v_1$};
\fill (1.9,3) node (v4) {} circle(3pt) ;
\node [below=2pt] at (1.95,3.03) {$v_4$};
\draw (v0)--node{$4$}(v1)--node{$10$}(v2)--node{$5$}(v0)--node{$8$}(v3)--node{$7$}(v2);
\draw (v0)--node{$9$}(v4)--node{$10$}(v3);
\fill (1.35,4.5) node (v5) {} circle(3pt) ;
\node [above=2pt] at (1.35,4.5) {$v_5$};
\fill (2.1,4.75) node (v6) {} circle(3pt) ;
\node [above=2pt] at (2.1,4.75) {$v_4^{(1)}$};
\draw (v3)--node{$6$}(v5)--node{$11$}(v4)--node{$3$}(v6)--node{$7$}(v5);
\fill (-2.5,2.4) node (v7) {} circle(3pt) ;
\node [left=2pt] at (-2.45,2.45) {$v_0^{(3)}$};
\draw (v0)--node{$3$}(v7)--node{$6$}(v1);
\fill (0.66,4.4) node (v8) {} circle(3pt) ;
\node [above=2pt] at (0.66,4.4) {$v_3^{(3)}$};
\draw (v3)--node{$5$}(v8)--node{$3$}(v5);
\fill (2.75,4.4) node (v9) {} circle(3pt) ;
\node [above=2pt] at (2.85,4.4) {$v_4^{(2)}$};
\draw (v4)--node{$4$}(v9)--node{$1$}(v6);
\fill (-2.95,1.65) node (v10) {} circle(3pt) ;
\node [left=2pt] at (-2.85,1.79) {$v_0^{(2)}$};
\draw (v0)--node{$2$}(v10)--node{$5$}(v7);
\fill (0.15,3.98) node (v11) {} circle(3pt) ;
\node [above=2pt] at (0.07,3.9) {$v_3^{(2)}$};
\draw (v3)--node{$4$}(v11)--node{$2$}(v8);
\fill (3.15,3.8) node (v12) {} circle(3pt) ;
\node [above=2pt] at (3.4,3.8) {$v_4^{(3)}$};
\draw (v4)--node{$5$}(v12)--node{$2$}(v9);
\fill (-3.3,1) node (v13) {} circle(3pt) ;
\node [left=2pt] at (-3.2,1.1) {$v_0^{(1)}$};
\draw (v0)--node{$1$}(v13)--node{$4$}(v10);
\fill (3.45,3.25) node (v14) {} circle(3pt) ;
\node [left=2pt] at (4.3,3.45) {$v_4^{(4)}$};
\draw (v4)--node{$6$}(v14)--node{$3$}(v12);
\fill (-0.15,3.4) node (v15) {} circle(3pt) ;
\node [left=2pt] at (0.04,3.63) {$v_3^{(1)}$};
\draw (v3)--node{$3$}(v15)--node{$1$}(v11);
\end{tikzpicture}
\end{minipage}}

\caption{Case 1: Star ~ edge ~coloring of $H_{n_\Delta}$ ($\Delta =4,5,6,7$) } \label{t1}
\end{figure}
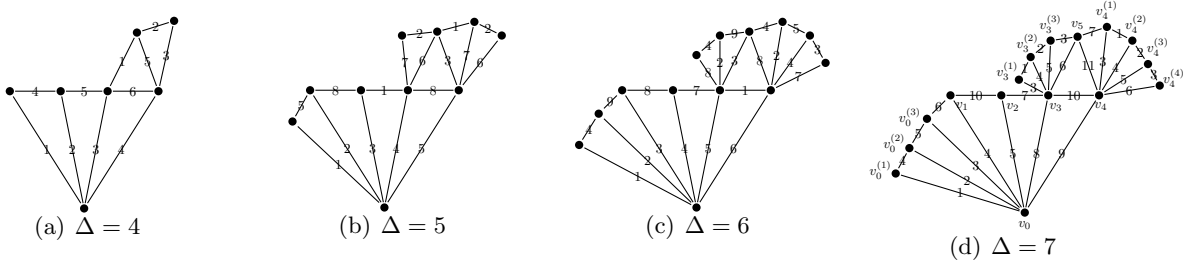
\indent Thus $\Delta \leq \chi_{st}^{'}(H_{n})  \leq \Delta+4.$ \\

\noindent{\bf Case~2:}\ \ This case of diagram is shown in Fig. 11,  we discuss the graph with $d(v_2)=d(v_3)=\Delta $ (denoted by $H^{2}_{n_\Delta}$). We consider the upper bound of $\chi_ {st}^{'}(H^{2}_{n_\Delta})$ based on the value of  $\Delta.$\\
 (1) When $\Delta =4,$ we consider the coloring of the Fig. 11(a), then $\chi_{st}^{'}(H^{2}_{n_\Delta}) \leq11= \Delta +3$.\\
 (2) If $\Delta =5,$ by the coloring of the Fig. 11(b), we know that $\chi_{st}^{'}(H^{2}_{n_\Delta}) \leq \Delta +3$, moreover $\chi_{st}^{'}(H_{n_\Delta}) \geq \chi_{st}^{'}(F_5)=\Delta +2,$ thus $\Delta +2=7\leq\chi_ {st}^{'}(H^{2}_{n_\Delta})\leq8=\Delta +3$.\\
 (3) When $\Delta =6,$ the Fig. 11(c) gives its a star edge coloring with $\Delta +3$ colors, thus $\chi_{st}^{'}(H^{2}_{n_\Delta}) \leq 9=\Delta +3$.\\
   (4) For $\Delta =7,$  considering the coloring of the Fig. 11(d), we obtain that $\chi_{st}^{'}(H^{2}_{n_\Delta}) \leq 10=\Delta +3$.\\
 (5) The case of $\Delta =8,$ a star edge coloring is shown in the Fig. 11(e), thus $\chi_{st}^{'}(H^{2}_{n_\Delta}) \leq11= \Delta +3$.\\
   (6) When $\Delta =9,$ we consider the coloring of the Fig. 11(f), then $\chi_{st}^{'}(H^{2}_{n_\Delta}) \leq11= \Delta +2$.\\
 (7) $\Delta =10,$ the Fig. 11(g) gives its a star edge coloring with $\Delta +1,$  thus $\chi_{st}^{'}(H^{2}_{n_\Delta}) \leq12= \Delta +2$.\\
(8) When $\Delta \geq 10,$ we consider the following coloring function $c:$ $c(v_0v_1)=\Delta+2$, $c(v_0v_2)=\Delta -1$, $c(v_0v_3)=\Delta +1$, $c(v_0v_4)=5$, $c(v_1v_2)=\Delta$, $c(v_2v_3)=\Delta -2$, $c(v_3v_4)=\Delta$, $c(v_1v_5)=\Delta +1$, $c(v_2v_5)=1$, $c(v_3v_6)=\Delta -3$, $c(v_4v_6)=\Delta +2$,  $c(v_2v_2^{(l)})=l+1$, $c(v_3v_3^{(m)})=m$, $c(v_5v_2^{(1)})=4$, $c(v_2^{(p)}v_2^{(p+1)})=p+4$, $c(v_2^{(\Delta -6)}v_2^{(\Delta -5)})=1$, $c(v_2^{(\Delta -5)}v_2^{(\Delta -4)})=2$, $c(v_3^{(q)}v_3^{(q+1)})=q+3$, $c(v_3^{(\Delta -5)}v_3^{(\Delta -4)})=1$, $c(v_3^{(\Delta -4)}v_4)=2$,  where $l,~m,~p,~q \in Z_{+}$ and $l,~ m \in [l,\Delta-4]$, $p,~q \in [l,\Delta-6]$.\\
\indent We prove the coloring function $c$ is a star edge coloring in the case(8): the paths $v_5-v_2^{(1)}-v_2^{(2)}-\cdots -v_2^{(\Delta -5)}-v_2^{(\Delta -4)}$, $v_3^{(1)}-v_3^{(2)}-\cdots -v_3^{(\Delta -5)}-v_3^{(\Delta -4)}-v_6$ are denoted by $P_2,$ $P_3$ respectively. By the theorem 1.4 of $^{[5]},$ we know that the Fans $P_2 \vee v_2$ and $P_3 \vee v_3$ have no bichromatic paths and cycles of length four. Assume the bichromatic path(cycle)$P$ containing an edge $v_0v_1$, $v_0v_2$, $v_0v_3$, $v_0v_4$, $v_1v_5$,$v_4v_6$,$v_1v_2$ or $v_3v_4$ in $H^{2}_{n_\Delta}.$ All the bichromatic 3-path(cycle) $P$ are $v_5-v_1-v_0-v_3$, $v_6-v_4-v_0-v_1$, $v_1-v_2-v_3-v_4$, $v_2^{(4)}-v_2-v_0-v_4$, $v_3^{(5)}-v_3-v_0-v_4$ in $H^{2}_{n_\Delta}.$  It is easy to check that any one of the above paths $P$ can't become bichromatic 4-path(cycle) by extending, then the graph $H^{2}_{n_\Delta}$ have no bichromatic 4-paths(cycles) including $v_0v_1$, $v_0v_2$, $v_0v_3$, $v_0v_4$, $v_1v_5$,$v_4v_6$,$v_1v_2$ or $v_3v_4.$  Since $c(v_2v_3)=\Delta-2$, and $\Delta -2$ appear only once in $H^{2}_{n_\Delta},$ thus $H^{2}_{n_\Delta}$ have no bichromatic 4-paths(cycles) including the edge $v_2v_3.$ By deleting the edges $v_0v_1$, $v_0v_2$, $v_0v_3$, $v_0v_4$, $v_1v_5$, $v_4v_6$,$v_1v_2$,$v_3v_4$ and $v_2v_3$, we obtain the subgraph $H^{2'}_{n_\Delta}$. The graph $H^{2'}_{n_\Delta}$ is not connected and has two connected components, which are Fans $P_2\vee v_2$ and $P_3\vee v_3$, thus $H^{2'}_{n_\Delta}$ have no bichromatic 4-paths(cycles). In summary,  $H^{2}_{n_\Delta}$ have no bichromatic 4-paths(cycles).  We obtain that $\chi_{st}^{'}(H^{2}_{n_\Delta}) \leq  \Delta+2$.

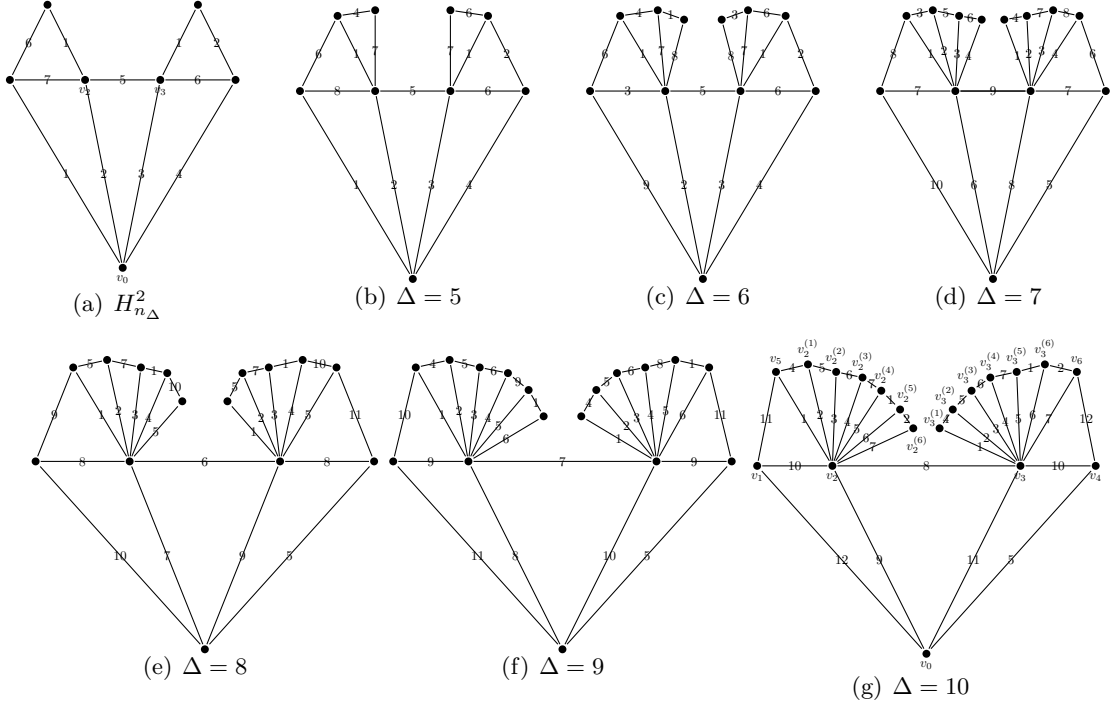
\begin{figure}[htpb]
\centering
\subfigure[$H^{2}_{n_\Delta}$]{
\begin{minipage}{3.6cm}
\centering
\begin{tikzpicture}[scale=0.5]
\tikzstyle{every node}=[font=\small,scale=0.5]
\fill (0,0) node (v0) {} circle(3pt) ;
\node [below=2pt] at (0,0) {$v_0$};
\fill (-1,5) node (v2) {} circle(3pt) ;
\node [below=2pt] at (-1,5) {$v_2$};
\fill (1,5) node (v3) {} circle(3pt) ;
\node [below=2pt] at (1,5) {$v_3$};
\fill (3,5) node (v4) {} circle(3pt) ;
\fill (-3,5) node (v1) {} circle(3pt) ;
\draw (v0)--node{$1$}(v1)--node{$7$}(v2)--node{$2$}(v0)--node{$3$}(v3)
--node{$5$}(v2);
\draw (v0)--node{$4$}(v4)--node{$6$}(v3);
\fill (-2,7) node (v5) {} circle(3pt) ;
\fill (2,7) node (v6) {} circle(3pt) ;
\draw (v1)--node{$6$}(v5)--node{$1$}(v2);
\draw (v3)--node{$1$}(v6)--node{$2$}(v4);
\end{tikzpicture}
\end{minipage}}
\subfigure[$\Delta =5$]{
\begin{minipage}{3.6cm}
\centering
\begin{tikzpicture}[scale=0.5]
\tikzstyle{every node}=[font=\small,scale=0.5]
\fill (0,0) node (v0) {} circle(3pt) ;
\fill (-1,5) node (v2) {} circle(3pt) ;
\fill (1,5) node (v3) {} circle(3pt) ;
\fill (3,5) node (v4) {} circle(3pt) ;
\fill (-3,5) node (v1) {} circle(3pt) ;
\draw (v0)--node{$1$}(v1)--node{$8$}(v2)--node{$2$}(v0)--node{$3$}(v3)
--node{$5$}(v2);
\draw (v0)--node{$4$}(v4)--node{$6$}(v3);
\fill (-2,7) node (v5) {} circle(3pt) ;
\fill (2,7) node (v6) {} circle(3pt) ;
\draw (v1)--node{$6$}(v5)--node{$1$}(v2);
\draw (v3)--node{$1$}(v6)--node{$2$}(v4);
\fill (-1,7.15) node (v7) {} circle(3pt) ;
\fill (1,7.15) node (v8) {} circle(3pt) ;
\draw (v5)--node{$4$}(v7)--node{$7$}(v2);
\draw (v6)--node{$6$}(v8)--node{$7$}(v3);
\end{tikzpicture}
\end{minipage}}
\subfigure[$\Delta =6$]{
\begin{minipage}{3.6cm}
\centering
\begin{tikzpicture}[scale=0.5]
\tikzstyle{every node}=[font=\small,scale=0.5]
\fill (0,0) node (v0) {} circle(3pt) ;
\fill (-1,5) node (v2) {} circle(3pt) ;
\fill (1,5) node (v3) {} circle(3pt) ;
\fill (3,5) node (v4) {} circle(3pt) ;
\fill (-3,5) node (v1) {} circle(3pt) ;
\draw (v0)--node{$9$}(v1)--node{$3$}(v2)--node{$2$}(v0)--node{$3$}(v3)
--node{$5$}(v2);
\draw (v0)--node{$4$}(v4)--node{$6$}(v3);
\fill (-2.2,7) node (v5) {} circle(3pt) ;
\fill (2.2,7) node (v6) {} circle(3pt) ;
\draw (v1)--node{$6$}(v5)--node{$1$}(v2);
\draw (v3)--node{$1$}(v6)--node{$2$}(v4);
\fill (-1.2,7.15) node (v7) {} circle(3pt) ;
\fill (1.2,7.15) node (v8) {} circle(3pt) ;
\draw (v5)--node{$4$}(v7)--node{$7$}(v2);
\draw (v6)--node{$6$}(v8)--node{$7$}(v3);
\fill (-0.5,6.9) node (v9) {} circle(3pt) ;
\fill (0.5,6.9) node (v10) {} circle(3pt) ;
\draw (v7)--node{$1$}(v9)--node{$8$}(v2);
\draw (v8)--node{$3$}(v10)--node{$8$}(v3);
\end{tikzpicture}
\end{minipage}}
\subfigure[$\Delta =7$]{
\begin{minipage}{3.6cm}
\centering
\begin{tikzpicture}[scale=0.5]
\tikzstyle{every node}=[font=\small,scale=0.5]

\fill (0,0) node (v0) {} circle(3pt) ;
\fill (-1,5) node (v2) {} circle(3pt) ;
\fill (1,5) node (v3) {} circle(3pt) ;
\fill (3,5) node (v4) {} circle(3pt) ;
\fill (-3,5) node (v1) {} circle(3pt) ;
\draw (v0)--node{$10$}(v1)--node{$7$}(v2)--node{$6$}(v0)--node{$8$}(v3)
--node{}(v2);
\draw  (v2)--node{$9$}(v3) ;
\draw (v0)--node{$5$}(v4)--node{$7$}(v3);
\fill (-2.3,7) node (v5) {} circle(3pt) ;
\fill (2.3,7) node (v6) {} circle(3pt) ;
\draw (v1)--node{$8$}(v5)--node{$1$}(v2);
\draw (v3)--node{$4$}(v6)--node{$6$}(v4);
\fill (-1.6,7.15) node (v7) {} circle(3pt) ;
\fill (1.6,7.15) node (v8) {} circle(3pt) ;
\draw (v5)--node{$3$}(v7)--node{$2$}(v2);
\draw (v6)--node{$8$}(v8)--node{$3$}(v3);
\fill (-0.9,7) node (v9) {} circle(3pt) ;
\fill (0.9,7) node (v10) {} circle(3pt) ;
\draw (v7)--node{$5$}(v9)--node{$3$}(v2);
\draw (v8)--node{$7$}(v10)--node{$2$}(v3);
\fill (-0.3,6.9) node (v11) {} circle(3pt) ;
\fill (0.3,6.9) node (v12) {} circle(3pt) ;
\draw (v9)--node{$6$}(v11)--node{$4$}(v2);
\draw (v10)--node{$4$}(v12)--node{$1$}(v3);
\end{tikzpicture}
\end{minipage}}
\subfigure[$\Delta =8$]{
\begin{minipage}{4.5cm}
\centering
\begin{tikzpicture}[scale=0.5]
\tikzstyle{every node}=[font=\small,scale=0.5]
\fill (0,0) node (v0) {} circle(3pt) ;
\fill (-2,5) node (v2) {} circle(3pt) ;
\fill (2,5) node (v3) {} circle(3pt) ;
\fill (4.5,5) node (v4) {} circle(3pt) ;
\fill (-4.5,5) node (v1) {} circle(3pt) ;
\draw (v0)--node{$10$}(v1)--node{$8$}(v2)--node{$7$}(v0)--node{$9$}(v3)
--node{$6$}(v2);
\draw (v0)--node{$5$}(v4)--node{$8$}(v3);
\fill (-3.5,7.5) node (v5) {} circle(3pt) ;
\fill (3.5,7.5) node (v6) {} circle(3pt) ;
\draw (v1)--node{$9$}(v5)--node{$1$}(v2);
\draw (v3)--node{$5$}(v6)--node{$11$}(v4);
\fill (-2.6,7.7) node (v7) {} circle(3pt) ;
\fill (2.6,7.7) node (v8) {} circle(3pt) ;
\draw (v5)--node{$5$}(v7)--node{$2$}(v2);
\draw (v6)--node{$10$}(v8)--node{$4$}(v3);
\fill (-1.7,7.5) node (v9) {} circle(3pt) ;
\fill (1.7,7.5) node (v10) {} circle(3pt) ;
\draw (v7)--node{$7$}(v9)--node{$3$}(v2);
\draw (v8)--node{$1$}(v10)--node{$3$}(v3);
\fill (-1,7.35) node (v11) {} circle(3pt) ;
\fill (1,7.35) node (v12) {} circle(3pt) ;
\draw (v9)--node{$1$}(v11)--node{$4$}(v2);
\draw (v10)--node{$7$}(v12)--node{$2$}(v3);
\fill (-0.6,6.6) node (v13) {} circle(3pt) ;
\fill (0.6,6.6) node (v14) {} circle(3pt) ;
\draw (v11)--node{$10$}(v13)--node{$5$}(v2);
\draw (v12)--node{$5$}(v14)--node{$1$}(v3);
\end{tikzpicture}
\end{minipage}}
\subfigure[$\Delta =9$]{
\begin{minipage}{4.5cm}
\centering
\begin{tikzpicture}[scale=0.5]
\tikzstyle{every node}=[font=\small,scale=0.5]
\fill (0,0) node (v0) {} circle(3pt) ;
\fill (-2.5,5) node (v2) {} circle(3pt) ;
\fill (2.5,5) node (v3) {} circle(3pt) ;
\fill (4.5,5) node (v4) {} circle(3pt) ;
\fill (-4.5,5) node (v1) {} circle(3pt) ;
\draw (v0)--node{$11$}(v1)--node{$9$}(v2)--node{$8$}(v0)--node{$10$}(v3)
--node{$7$}(v2);
\draw (v0)--node{$5$}(v4)--node{$9$}(v3);
\fill (-3.9,7.5) node (v5) {} circle(3pt) ;
\fill (3.9,7.5) node (v6) {} circle(3pt) ;
\draw (v1)--node{$10$}(v5)--node{$1$}(v2);
\draw (v3)--node{$6$}(v6)--node{$11$}(v4);
\fill (-3,7.7) node (v7) {} circle(3pt) ;
\fill (3,7.7) node (v8) {} circle(3pt) ;
\draw (v5)--node{$4$}(v7)--node{$2$}(v2);
\draw (v6)--node{$1$}(v8)--node{$5$}(v3);
\fill (-2.2,7.5) node (v9) {} circle(3pt) ;
\fill (2.2,7.5) node (v10) {} circle(3pt) ;
\draw (v7)--node{$5$}(v9)--node{$3$}(v2);
\draw (v8)--node{$8$}(v10)--node{$4$}(v3);
\fill (-1.45,7.35) node (v11) {} circle(3pt) ;
\fill (1.45,7.35) node (v12) {} circle(3pt) ;
\draw (v9)--node{$6$}(v11)--node{$4$}(v2);
\draw (v10)--node{$6$}(v12)--node{$3$}(v3);
\fill (-0.9,6.9) node (v13) {} circle(3pt) ;
\fill (0.9,6.9) node (v14) {} circle(3pt) ;
\draw (v11)--node{$9$}(v13)--node{$5$}(v2);
\draw (v12)--node{$5$}(v14)--node{$2$}(v3);
\fill (-0.5,6.2) node (v15) {} circle(3pt) ;
\fill (0.5,6.2) node (v16) {} circle(3pt) ;
\draw (v13)--node{$1$}(v15)--node{$6$}(v2);
\draw (v14)--node{$4$}(v16)--node{$1$}(v3);
\end{tikzpicture}
\end{minipage}}
\subfigure[$\Delta =10$]{
\begin{minipage}{4.5cm}
\centering
\begin{tikzpicture}[scale=0.5]
\tikzstyle{every node}=[font=\small,scale=0.5]
\fill (0,0) node (v0) {} circle(3pt) ;
\node [below=2pt] at (0,0) {$v_0$};
\fill (-2.5,5) node (v2) {} circle(3pt) ;
\node [below=2pt] at(-2.5,5) {$v_2$};
\fill (2.5,5) node (v3) {} circle(3pt) ;
\node [below=2pt] at(2.5,5) {$v_3$};
\fill (4.5,5) node (v4) {} circle(3pt) ;
\node [below=2pt] at(4.5,5) {$v_4$};
\fill (-4.5,5) node (v1) {} circle(3pt) ;
\node [below=2pt] at(-4.5,5) {$v_1$};
\draw (v0)--node{$12$}(v1)--node{$10$}(v2)--node{$9$}(v0)--node{$11$}(v3)
--node{$8$}(v2);
\draw (v0)--node{$5$}(v4)--node{$10$}(v3);
\fill (-4,7.5) node (v5) {} circle(3pt) ;
\node [above=2pt] at(-4,7.5) {$v_5$};
\fill (4,7.5) node (v6) {} circle(3pt) ;
\node [above=2pt] at(4,7.5) {$v_6$};
\draw (v1)--node{$11$}(v5)--node{$1$}(v2);
\draw (v3)--node{$7$}(v6)--node{$12$}(v4);
\fill (-3.15,7.7) node (v7) {} circle(3pt) ;
\node [above=2pt] at(-3.15,7.7) {$v_2^{(1)}$};
\fill (3.15,7.7) node (v8) {} circle(3pt) ;
\node [above=2pt] at(3.15,7.7) {$v_3^{(6)}$};
\draw (v5)--node{$4$}(v7)--node{$2$}(v2);
\draw (v6)--node{$2$}(v8)--node{$6$}(v3);
\fill (-2.4,7.5) node (v9) {} circle(3pt) ;
\node [above=2pt] at(-2.4,7.5) {$v_2^{(2)}$};
\fill (2.4,7.5) node (v10) {} circle(3pt) ;
\node [above=2pt] at(2.4,7.5) {$v_3^{(5)}$};
\draw (v7)--node{$5$}(v9)--node{$3$}(v2);
\draw (v8)--node{$1$}(v10)--node{$5$}(v3);
\fill (-1.7,7.35) node (v11) {} circle(3pt) ;
\node [above=2pt] at(-1.7,7.35){$v_2^{(3)}$};
\fill (1.7,7.35) node (v12) {} circle(3pt) ;
\node [above=2pt] at(1.7,7.35){$v_3^{(4)}$};
\draw (v9)--node{$6$}(v11)--node{$4$}(v2);
\draw (v10)--node{$7$}(v12)--node{$4$}(v3);
\fill (-1.2,7) node (v13) {} circle(3pt) ;
\node [above=2pt] at(-1.1,7) {$v_2^{(4)}$};
\fill (1.2,7) node (v14) {} circle(3pt) ;
\node [above=2pt] at(1.1,7) {$v_3^{(3)}$};
\draw (v11)--node{$7$}(v13)--node{$5$}(v2);
\draw (v12)--node{$6$}(v14)--node{$3$}(v3);
\fill (-0.7,6.5) node (v15) {} circle(3pt) ;
\node [above=2pt] at(-0.5,6.5) {$v_2^{(5)}$};
\fill (0.7,6.5) node (v16) {} circle(3pt) ;
\node [above=2pt] at(0.5,6.45) {$v_3^{(2)}$};
\draw (v13)--node{$1$}(v15)--node{$6$}(v2);
\draw (v14)--node{$5$}(v16)--node{$2$}(v3);
\fill (-0.35,6) node (v17) {} circle(3pt) ;
\node [below=2pt] at(-0.25,6) {$v_2^{(6)}$};
\fill (0.35,6) node (v18) {} circle(3pt) ;
\node [above=2pt] at(0.2,5.85) {$v_3^{(1)}$};
\draw (v15)--node{$2$}(v17)--node{$7$}(v2);
\draw (v16)--node{$4$}(v18)--node{$1$}(v3);
\end{tikzpicture}
\end{minipage}}
\caption{Case 2: the star edge coloring of $H^{2}_{n_\Delta}$($\Delta =5,\cdots ,10$)} \label{t1}
\end{figure}

\indent Thus, $\Delta \leq \chi_{st}^{'}(H_{n}) \leq \chi_{st}^{'}(H^{2}_{n_\Delta}) \leq \Delta+2$. \\
\indent According to the Cases 3.5.1 and 3.5.2, when $H_{n}\in\xi_n^{3},$ we obtain that  $\Delta(H_{n})\leq \chi_{st}^{'}(H_{n})\leq\Delta(H_{n})+4.$ In the document\cite{DS}, a classification of $G_{n}(\in\zeta_n^{3})$ has been given. On the basis of the classification, we find that for every $G_{n}\in\zeta_n^{3}$ there exists a $H_{n} \in \xi_n^{3}$ such that $G_{n}$ is a subgraph of $H_{n}$ and $\Delta(H_{n})\leq\Delta(G_n)+2.$ Therefore, we can obtain that  $\Delta(G_{n})\leq\chi_{st}^{'}(G_{n}) \leq \chi_{st}^{'}(H_{n})\leq\Delta(H_{n})+4\leq\Delta(G_{n})+6.$ $\hfill \Box $\\
\section{Outerplanar graphs with $\Delta =5$ }

\noindent{\bf Theorem $4.1$} Let $G$ be a 2-connected outerplanar graph with $n$ vertices and $\Delta(G)=5,$ then $\chi_{st}^{'}(G)\leq\lfloor 1.5\Delta \rfloor + 2=\Delta+4=9.$\\

\noindent{\bf Proof.} Since $G$ has $n$ vertices and $\Delta(G)=5,$ obviously $n\geq 6.$ When $n=6,$ $G$ is isomorphic to $F_6.$ By Lemma 2.3, we know that $\chi_{st}^{'}(F_6)=6<9.$
Since maximal outerplanar graphs with maximal degree 5 can be drawn as Fig. 12 when the order is going infinite and the $F_6s$ between $E$ and $F$ form a cycle. Inspired by the document\cite{YS}, we give a star 9-edge coloring as Fig. 12 showing. If $n>6,$ any outerplanar graph with maximum degree 5 can become a maximal outerplanar graph with maximum degree 5 by adding edges, moreover it is a subgraph of Fig. 12. Since Fig. 12 has a 9 star edge coloring, we obtain that any 2-connected outerplanar graph with $\Delta(G)=5$ has star chromatic index at most 9.
 \vskip 0.2cm
\begin{figure}[htpb]
\centering
\usetikzlibrary{decorations.pathreplacing}
\begin{tikzpicture}[scale=0.8]
\tikzstyle{every node}=[font=\small,scale=0.8]
\node (v5) at (-2,0) {};
\node (v4) at (-3,0) {};
\node (v3) at (-4,0) {};
\node (v1) at (-5,0) {};
\node (v2) at (-5,1) {};
\node (v6) at (-3.5,1.5) {};
\draw  (v1)  --node [pos=0.5,sloped]{$2$} (v2);
\draw  (v1)  --node [pos=0.5,sloped]{$5$} (v3);
\draw  (v3)  --node [pos=0.5,sloped]{$a$} (v4);
\draw  (v4)  --node [pos=0.5,sloped]{$6$}(v5);
\draw  (v5) --node [pos=0.5,sloped]{$a$} (v6);
\draw  (v6)  --node [pos=0.5,sloped]{$c$} (v2);
\draw  (v6)  --node [pos=0.5,sloped]{$4$} (v1);
\draw  (v6)  --node [pos=0.5,sloped]{$3$} (v3);
\draw  (v6)  --node [pos=0.5,sloped]{$1$} (v4);
\node (v9) at (-5,-1.5) {};
\node (v8) at (-6,0) {};
\node (v7) at (-7,0) {};
\node (v10) at (-3.5,-1) {};
\draw  (v7)  --node [pos=0.5,sloped]{$5$} (v8);
\draw  (v8)  --node [pos=0.5,sloped]{$b$} (v1);
\draw  (v1)  --node [pos=0.5,sloped]{$6$} (v9);
\draw  (v9)  --node [pos=0.5,sloped]{$b$} (v10);
\draw  (v10)  --node [pos=0.5,sloped]{$4$} (v3);
\draw  (v3)  --node [pos=0.5,sloped]{$c$} (v9);
\draw  (v9)  --node [pos=0.5,sloped]{$a$} (v8);
\draw  (v9)  --node [pos=0.5,sloped]{$1$} (v7);
\node (v12) at (-2,-1.5) {};
\node (v11) at (-3,-1) {};
\node (v13) at (-1,0) {};
\node (v14) at (-0.5,-1) {};
\draw  (v4)  --node [pos=0.5,sloped]{$b$} (v11);
\draw  (v11)  --node [pos=0.5,sloped]{$5$} (v12);
\draw  (v12)  --node [pos=0.5,sloped]{$4$} (v4);
\draw  (v12)  --node [pos=0.5,sloped]{$c$} (v5);
\draw  (v5)  --node [pos=0.5,sloped]{$b$} (v13);
\draw  (v13)  --node [pos=0.5,sloped]{$5$} (v14);
\draw  (v14) --node [pos=0.5,sloped]{$6$}(v12);
\draw  (v12)  --node [pos=0.5,sloped]{$3$} (v13);
\node (v15) at (-1,1.5) {};
\node (v16) at (0,0) {};
\node (v17) at (1,0) {};
\node (v18) at (1,1) {};
\draw  (v15)  --node [pos=0.5,sloped]{$2$} (v5);
\draw  (v15)  --node [pos=0.5,sloped]{$4$} (v13);
\draw  (v13)  --node [pos=0.5,sloped]{$a$} (v16);
\draw  (v16)  --node [pos=0.5,sloped]{$1$} (v17);
\draw  (v17)  --node [pos=0.5,sloped]{$3$} (v18);
\draw  (v15)  --node [pos=0.5,sloped]{$5$} (v18);
\draw  (v15)  --node [pos=0.5,sloped]{$b$} (v17);
\draw  (v15)  --node [pos=0.5,sloped]{$c$}(v16);
\node (v20) at (1,-1.5) {};
\node (v19) at (0,-1) {};
\node (v21) at (2,0) {};
\node (v22) at (3,0) {};
\draw (v20) --node [pos=0.5,sloped]{$c$} (v19);
\draw  (v19) --node [pos=0.5,sloped]{$3$} (v16);
\draw  (v16)  --node [pos=0.5,sloped]{$6$} (v20);
\draw  (v20)  --node [pos=0.5,sloped]{$5$} (v17);
\draw  (v17)  --node [pos=0.5,sloped]{$a$} (v21);
\draw  (v21) --node [pos=0.5,sloped]{$c$} (v22);
\draw  (v22)  --node [pos=0.5,sloped]{$2$} (v20);
\draw  (v20)  --node [pos=0.5,sloped]{$4$} (v21);
\node (v23) at (3,1.5) {};
\node (v24) at (1.5,1) {};
\node (v25) at (4,0) {};
\draw  (v23)  --node [pos=0.5,sloped]{$6$}(v24);
\draw  (v24)  --node [pos=0.5,sloped]{$5$} (v21);
\draw  (v21)  --node [pos=0.5,sloped]{$b$}(v23);
\draw  (v22)  --node [pos=0.5,sloped]{$3$} (v23);
\draw  (v22)  --node [pos=0.5,sloped]{$5$} (v25);
\draw  (v23)  --node [pos=0.5,sloped]{$1$} (v25);
\node (v27) at (-7,1.5) {};
\node (v29) at (-8,0) {};
\node (v28) at (-8,1) {};
\node (v30) at (-5.5,1) {};
\draw  (v27)  --node [pos=0.5,sloped]{$6$} (v28);
\draw  (v28)  --node [pos=0.5,sloped]{$5$} (v29);
\draw  (v29)  --node [pos=0.5,sloped]{$c$} (v7);
\draw  (v30)  --node [pos=0.5,sloped]{$4$} (v8);
\draw  (v27)  --node [pos=0.5,sloped]{$b$} (v29);
\draw  (v27)  --node [pos=0.5,sloped]{$3$} (v7);
\draw  (v27)  --node [pos=0.5,sloped]{$1$} (v8);
\draw  (v27)  --node [pos=0.5,sloped]{$c$}(v30);
\node (v31) at (-8,-1.5) {};
\node (v32) at (-9,0) {};
\node (v34) at (-10,0) {};
\node (v33) at (-10,-1) {};
\draw (v31) --node [pos=0.5,sloped]{$2$} (v7);
\draw (v31) --node [pos=0.5,sloped]{$4$} (v29);
\draw  (v29) --node [pos=0.5,sloped]{$a$} (v32);
\draw (v31) --node [pos=0.5,sloped]{$5$} (v32);
\draw (v31) --node [pos=0.5,sloped]{$6$} (v34);
\draw (v31) --node [pos=0.5,sloped]{$c$} (v33);
\draw  (v33) edge (v34);
\draw (v34) --node [pos=0.5,sloped]{$3$} (v33);
\draw  (v34)  --node [pos=0.5,sloped]{$1$} (v32);
\node (v35) at (4.5,1) {};
\draw  (v35)  --node [pos=0.5,sloped]{$4$} (v25);
\draw (v23)  --node [pos=0.5,sloped]{$c$} (v35);
\node (v36) at (5,0) {};
\node (v37) at (6,0) {};
\node (v38) at (6.5,-1) {};
\node (v26) at (4,-1.5) {};
\draw (v22)  --node [pos=0.5,sloped]{$1$} (v26);
\draw(v26)  --node [pos=0.5,sloped]{$a$} (v25);
\draw(v25)  --node [pos=0.5,sloped]{$b$} (v36);
\draw(v36)  --node [pos=0.5,sloped]{$5$} (v37);
\draw (v37)  --node [pos=0.5,sloped]{$4$} (v38);
\draw(v38)  --node [pos=0.5,sloped]{$b$} (v26);
\draw (v26)  --node [pos=0.5,sloped]{$c$} (v37);
\draw (v36)  --node [pos=0.5,sloped]{$6$} (v26);
\node (v41) at (-11,0) {};
\node (v39) at (-12.5,0) {};
\node (v42) at (-9,1) {};
\node (v40) at (-10,1.5) {};
\draw  (v39)  --node [pos=0.5,sloped]{$2$} (v40);
\draw (v40)  --node [pos=0.5,sloped]{$4$} (v41);
\draw (v41)  --node [pos=0.5,sloped]{$b$} (v39);
\draw (v41)  --node [pos=0.5,sloped]{$a$} (v34);
\draw (v40)  --node [pos=0.5,sloped]{$c$} (v34);
\draw (v40)  --node [pos=0.5,sloped]{$b$} (v32);
\draw  (v42)  --node [pos=0.5,sloped]{$3$} (v32);
\draw (v42)  --node [pos=0.5,sloped]{$5$} (v40);
\draw[ decorate, decoration={brace, raise=5pt,amplitude=5pt}]  (-7,1.5) --node[pos=0.5,above=9,sloped]{$E$}(3,1.5);
\draw[decorate,decoration={brace,mirror,raise=5pt}] (-8,-1.5) -- node[pos=0.5,below=6,sloped]{$F$}(1,-1.5);
\foreach \p in {v1,v2,v3,v4,v5,v6,v7,v8,v9,v10,
v11,v12,v13,v14,v15,v16,v17,v18,v19,v20,v21,
v22,v23,v24,v25,v26,v27,v28,v29,v30,
v31,v32,v33,v34,v35,v36,v37,v38,v39,v40,v41,v42} \fill
[opacity=0.75](\p) circle(3pt);
\end{tikzpicture}
\caption{$\mathrm{A~9~star~edge~coloring}$ }\label{fig1}
\end{figure}
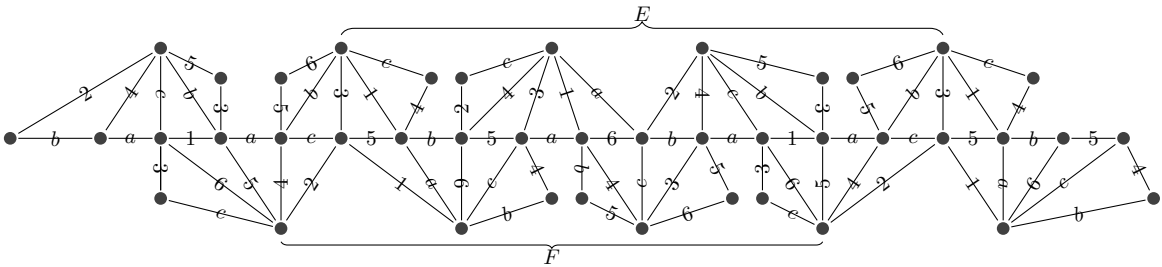

\begin{figure}[htpb]
\centering

\end{figure}

\section{Concluding remarks}

 \vspace{2mm}
In this paper, we discuss the star chromatic index of outerplanar graphs with small diameter inspired by the Conjecture 1.5. Otherwise, we find that the upper bound of star chromatic index of outerplanar graphs with diameter 2 or 3 is
$\Delta+6,$ which is much less than $\lfloor\frac{3\Delta}{2}\rfloor+1$ when $\Delta$ increases infinitely. Theorem 1.5 tells us that $\chi^{\prime}_{st}(G)\leq 5$ when $G$ is a subcubic outerplanar graph.  By Theorem 1.7, we know that outerplanar  graphs with $\Delta=4$ have star chromatic index 6. Moreover, Theorem 4.1 obtains that 2-connected outerplanar graphs have star chromatic indexes no more than 9 when $\Delta=5.$  So overall, these results make us sure that the following conjecture is true.\\

\noindent{\bf Conjecture 4.1}
If $G$ is a 2-connected outerplanar graph with the maximum degree $\Delta\geq 6$, then $\chi_{st}^{'} (G) \leq \Delta  + 6.$\\

\noindent{\bf Conjecture 4.2}
If $G$ is a 2-connected maximal outerplanar graph with the maximum degree $\Delta\geq 6$, then $\chi_{st}^{'} (G) \leq \Delta  + 4.$\\

In \cite{LB}, Bezegov$\mathrm{\acute{a}}$ et al. proved that $\chi_{st}^{'} (T) \leq \lfloor\frac{3\Delta}{2}\rfloor,$ when T be a tree with maximum degree $\Delta$ and every tree with a $\Delta$-vertex whose all neighbors are $\Delta$-vertices achieves the upper bound. We propose the following open question inspired by Bezegov$\mathrm{\acute{a}}$'s and our results.\\

\noindent{\bf Problem 4.2}
If $G$ is a 2-connected graph with linear size and  maximum degree $\Delta$, then $\chi_{st}^{'} (G) \leq \Delta  + c$ (c is a constant)?\\

\end{document}